# Pauli Pascal Pyramids, Pauli Fibonacci Numbers, and Pauli Jacobsthal Numbers

Martin Erik Horn*

E-Mail: *horn@fhw-berlin.de*
*mail@grassmann-algebra.de*

**Abstract**

The three anti-commutative two-dimensional Pauli Pascal triangles can be generalized into multi-dimensional Pauli Pascal hyperpyramids.

Fibonacci and Jacobsthal numbers are then generalized into Pauli Fibonacci numbers, Pauli Jacobsthal numbers, and Pauli Fibonacci numbers of higher order.

And the question is: are Pauli rabbits killer rabbits?

**Contents**

* Address: FHW Berlin, Badensche Str. 50/51, Fach # 99, D – 10825 Berlin,
Otto-Hahn-Schule Neukölln, Buschkrugallee 63, D – 12359 Berlin.

## 1. Introduction

Mathematics is an astonishing subject. You never know, what the basis of mathematics is. Is mathematics logical? Is mathematics fallible? And yet, there are mathematical constructions, which are fascinatingly beautiful and deeply impressive. They surely reflect something of a hidden truth.

As physicist and physics education specialist I look at these mathematical relations – like for example Ramanujans ${}_1\psi_1$–summation formula – with amazement and great astonishment and ask myself: This formula can be found? How can a formula like that be found? Finding new mathematical relations is certainly far more difficult than to prove them after the process of discovery or invention.

But the density of mathematical texts often prevents readers from understanding the actual, deeper meaning of these texts. They are able to understand the mathematical letters and symbols. Yet these letters and symbols are combined into a narrow and one-dimensional mathematical language. So they understand and can follow the first layer of the text and grasp the mathematical statements only superficially. But there is something behind mathematical texts, a deeper meaning and deeper context, which very often is not accessible to the readers.

On my way through the mysterious and hidden worlds of mathematics I therefore prefer to travel like John S. Bell, who once remarked that "the longer road sometimes gives more familiarity with the country" [10]. The surprising landscape which can be discovered behind unexpected bends of the road then may help to understand, how something can be discovered and how modern mathematics may be deciphered. In the following I present first steps on the very long road through a fascinating country: the country of Pascal space and Pauli Pascal space .

Of course lots of people already entered this Pascalian landscape. They found and described the Pascal triangle, the Pascal tetrahedron, the Pascal pyramid and Pascal simplexes, called Pascal hyperpyramids in this text. And they tell us about interesting relations between these Pascal structures and Fibonacci numbers or Jacobsthal numbers (for example see [1] ). But – in two- and three-dimensional situations – you can draw these triangles and pyramids. Therefore the landscape opens in front of you. And then you can indeed see mathematical relations which were hidden in the dark world of dense mathematical texts without illustrations before. This is just done here: illustrations of Pascal space are given for our world of ordinary numbers and for the world of Pauli numbers.

## 2. Pascal Triangles

The Taylor expansion (Mac Laurin form) of $(a + b)^n$ clearly shows that there are three Pascal triangles.

For $n \geq 0$ there is:

$$\begin{aligned}
(a + b)^0 &= 1 \\
(a + b)^1 &= 1a + 1b \\
(a + b)^2 &= 1a^2 + 2ab + 1b^2 \\
(a + b)^3 &= 1a^3 + 3a^2b + 3ab^2 + 1b^3 \\
(a + b)^4 &= 1a^4 + 4a^3b + 6a^2b^2 + 4ab^3 + 1b^4 \\
&\text{etc…}
\end{aligned} \tag{1}$$

For $n < 0$ and $|a| > |b|$ there is:

$$
\begin{aligned}
(a + b)^{-1} &= 1a^{-1} - 1a^{-2}b + 1a^{-3}b^{2} - la^{-4}b^{3} + \ldots - \ldots \\
(a + b)^{-2} &= 1a^{-2} - 2a^{-3}b + 3a^{-4}b^{2} - 4a^{-5}b^{3} + \ldots - \ldots \\
(a + b)^{-3} &= 1a^{-3} - 3a^{-4}b + 6a^{-5}b^{2} - 10a^{-6}b^{3} + \ldots - \ldots \\
(a + b)^{-4} &= 1a^{-4} - 4a^{-5}b + 10a^{-6}b^{2} - 20a^{-7}b^{3} + \ldots - \ldots \\
(a + b)^{-5} &= 1a^{-5} - 5a^{-6}b + 15a^{-7}b^{2} - 35a^{-8}b^{3} + \ldots - \ldots \\
&\text{etc}\ldots
\end{aligned}
\qquad (2)
$$

For $n < 0$ and $|a| < |b|$ there is:

$$
\begin{aligned}
(a + b)^{-1} &= 1a^{0}b^{-1} - la^{1}b^{-2} + 1a^{2}b^{-3} - la^{3}b^{-4} + \ldots - \ldots \\
(a + b)^{-2} &= 1a^{0}b^{-2} - 2a^{1}b^{-3} + 3a^{2}b^{-4} - 4a^{3}b^{-5} + \ldots - \ldots \\
(a + b)^{-3} &= 1a^{0}b^{-3} - 3a^{1}b^{-4} + 6a^{2}b^{-5} - l0a^{3}b^{-6} + \ldots - \ldots \\
(a + b)^{-4} &= 1a^{0}b^{-4} - 4a^{1}b^{-5} + 10a^{2}b^{-6} - 20a^{3}b^{-7} + \ldots - \ldots \\
(a + b)^{-5} &= 1a^{0}b^{-5} - 5a^{1}b^{-6} + 15a^{2}b^{-7} - 35a^{3}b^{-8} + \ldots - \ldots \\
&\text{etc}\ldots
\end{aligned}
\qquad (3)
$$

These coefficients of $(a + b)^{x+y}$ can be arranged symmetrically by describing exponents of $a$ as x-coordinate and exponents of $b$ as y-coordinate. If positive directions of this coordinate system are directed downwards the structure of the Pascal Triangles will be seen clearly (see figure 1). All coefficients of constant $n = x + y$ are then situated in horizontal lines.

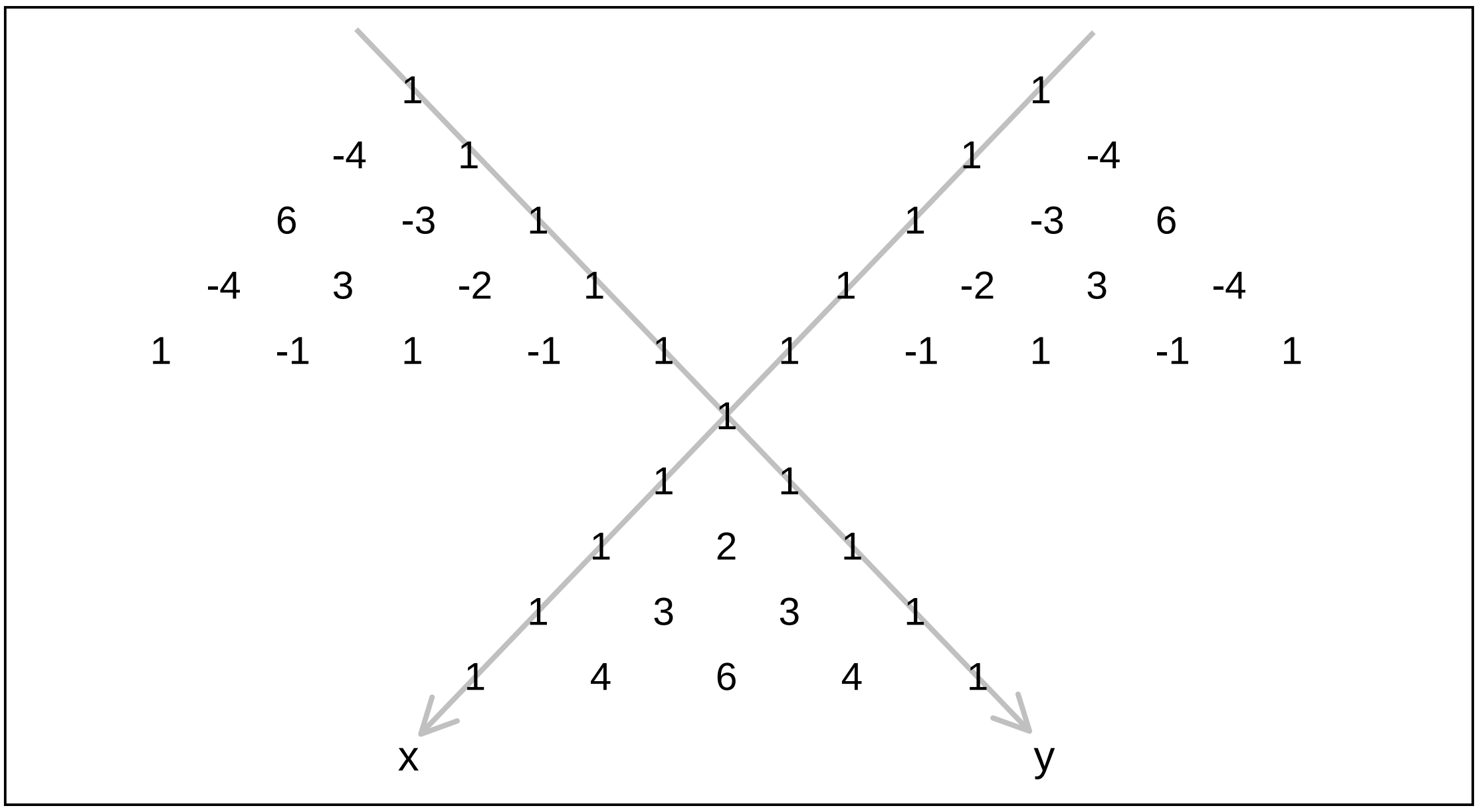

Fig. 1: The Pascal plane with the Pascal triangles.

More about the relation between these numbers of Pascal triangles and bilateral hypergeometric series can be found in [6] and [2].

### 3. Pascal Pyramids and Pascal Hyperpyramids

The Taylor expansion of $(a + b + c)^n$ gives the coefficients of four Pascal pyramids.

For $n \geq 0$ there is:

$$
\begin{aligned}
(a + b + c)^0 &= 1 \\
(a + b + c)^1 &= 1a + 1b + 1c \\
(a + b + c)^2 &= 1a^2 + 2ab + 1b^2 + 2bc + 1c^2 + 2ca \\
(a + b + c)^3 &= 1a^3 + 3a^2b + 3ab^2 + 1b^3 + 3b^2c + 3bc^2 + 1c^3 + 3c^2a + 3ca^2 + 6abc \\
(a + b + c)^4 &= 1a^4 + 4a^3b + 6a^2b^2 + 4ab^3 + 1b^4 + 4b^3c + 6b^2c^2 + 4bc^3 + 1c^4 \\
&\quad + 4c^3a + 6c^2a^2 + 4ca^3 + 12a^2bc + 12ab^2c + 12abc^2 \\
&\text{etc...}
\end{aligned}
\tag{4}
$$

These coefficients of $(a + b)^{x+y+z}$ can be arranged symmetrically by choosing exponents of $a$ as x-coordinate, exponents of $b$ as y-coordinate, and exponents of $c$ as z-coordinate. The three-dimensional Pascal pyramid in the positive region of Pascal space is shown in figure 2. The surface of this pyramid consists of three different Pascal triangles which lie in the xy-plane, the yz-plane, and the zx-plane.

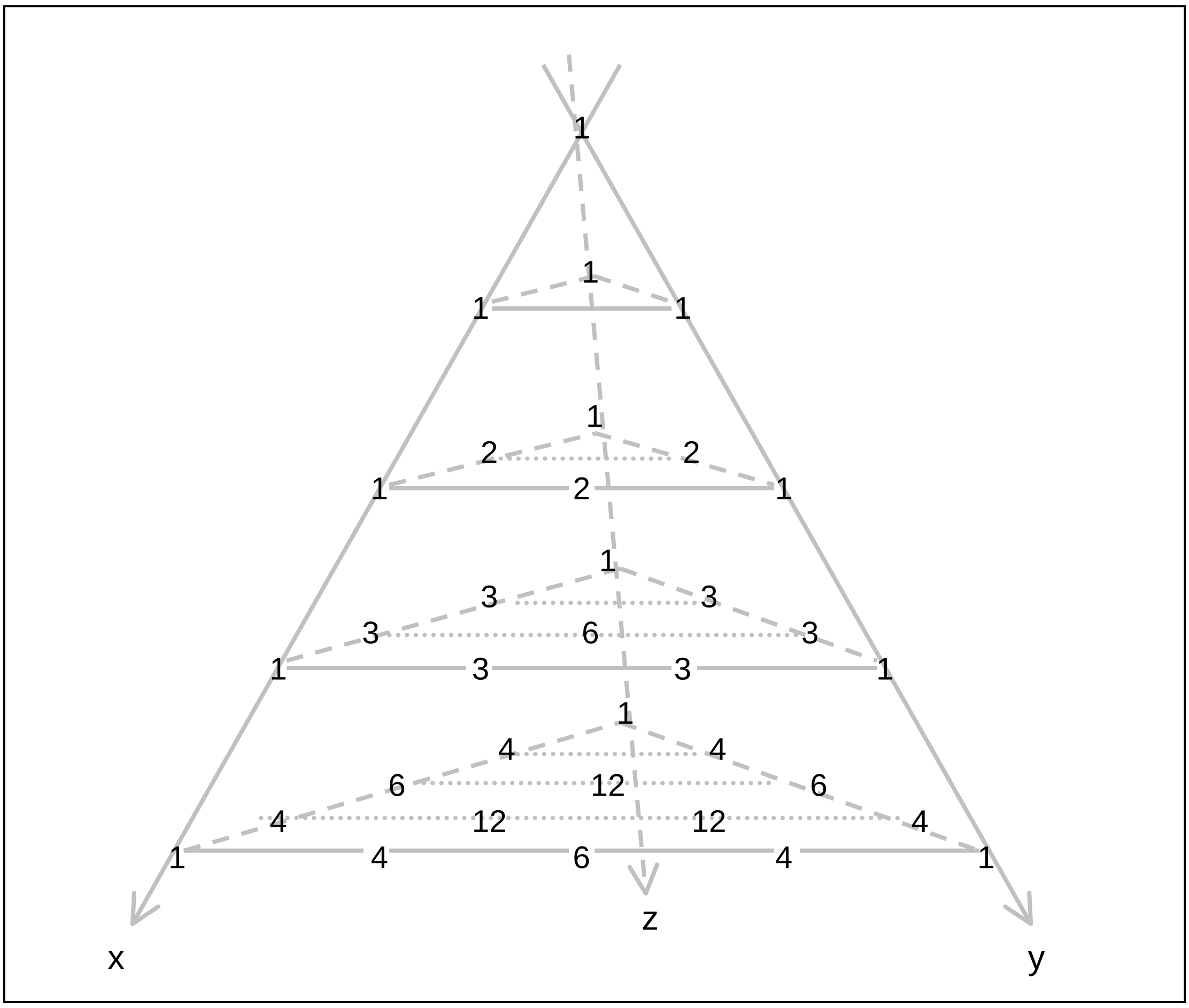

Fig. 2: Pascal space with the positive Pascal pyramid.

For $n < 0$ and $|a| > |b|, |c|$ there is:

$$\begin{aligned}
(a+b+c)^{-1} &= 1a^{-1} - 1a^{-2}b - 1a^{-2}c + 1a^{-3}b^2 + 2a^{-3}bc + 1a^{-3}c^2 \\
&\quad - 1a^{-4}b^3 - 3a^{-4}b^2c - 3a^{-4}bc^2 - 1a^{-4}c^3 + \ldots - \ldots \\
(a+b+c)^{-2} &= 1a^{-2} - 2a^{-3}b - 2a^{-3}c + 3a^{-4}b^2 + 6a^{-4}bc + 3a^{-4}c^2 \\
&\quad - 4a^{-5}b^3 - 12a^{-5}b^2c - 12a^{-5}bc^2 - 4a^{-5}c^3 + \ldots - \ldots \\
(a+b+c)^{-3} &= 1a^{-3} - 3a^{-4}b - 3a^{-4}c + 6a^{-5}b^2 + 12a^{-5}bc + 6a^{-5}c^2 \\
&\quad - 10a^{-6}b^3 - 30a^{-6}b^2c - 30a^{-6}bc^2 - 10a^{-6}c^3 + \ldots - \ldots \\
(a+b+c)^{-4} &= 1a^{-4} - 4a^{-5}b - 4a^{-5}c + 10a^{-6}b^2 + 20a^{-6}bc + 10a^{-6}c^2 \\
&\quad - 20a^{-7}b^3 - 60a^{-7}b^2c - 60a^{-7}bc^2 - 20a^{-7}c^3 + \ldots - \ldots \\
(a+b+c)^{-5} &= 1a^{-5} - 5a^{-6}b - 5a^{-6}c + 15a^{-7}b^2 + 30a^{-7}bc + 15a^{-7}c^2 \\
&\quad - 35a^{-8}b^3 - 105a^{-8}b^2c - 105a^{-8}bc^2 - 35a^{-8}c^3 + \ldots - \ldots \\
&\text{etc}\ldots
\end{aligned} \tag{5}$$

Thus the first negative Pascal pyramid occupies the region of the negative x-coordinate, positive y-coordinate, and positive z-coordinate (see figure 3).

The other two negative Pascal pyramids can be constructed in a similar way for $n < 0$. Then $|b| > |a|, |c|$ gives the second negative Pascal pyramid. This pyramid occupies the region of the negative y-coordinate, positive z-coordinate, and positive x-coordinate (see figure 4).

The third negative Pascal pyramid with $|c| > |a|, |b|$ occupies the region of the negative z-coordinate, positive x-coordinate, and positive y-coordinate, which is indicated in figure 5.

While the surface of the positive Pascal pyramid is composed of the positive Pascal triangle, the surfaces of negative Pascal pyramids are composed of negative Pascal triangles.

In four-dimensional Pascal hyperspace there exist one positive four-dimensional Pascal hyperpyramid and four negative four-dimensional Pascal hyperpyramids. The positive Pascal hyperpyramid can be constructed with the coefficients of the Taylor expansion of

$$(a + b + c + d)^n \quad \text{with } n \geq 0$$

or alternatively with the multinomial coefficients [8]. The four three-dimensional hypersurfaces of this positive four-dimensional Pascal hyperpyramid are the positive Pascal pyramid of figure 2.

The negative Pascal hyperpyramids can be constructed with the coefficients of the Taylor expansion of $(a + b + c + d)^n$ with $n < 0$. The four three-dimensional hypersurfaces of these negative four-dimensional Pascal hyperpyramids are the negative Pascal pyramids of figure 2.

In a similar way $(k + 1)$ Pascal hyperpyramids of dimension k can be constructed with the coefficients of the Taylor expansion of $(a_1 + a_2 + a_3 + \ldots + a_k)^n$. Unfortunately these hyperpyramids can not be drawn in our 3-dimensional world.

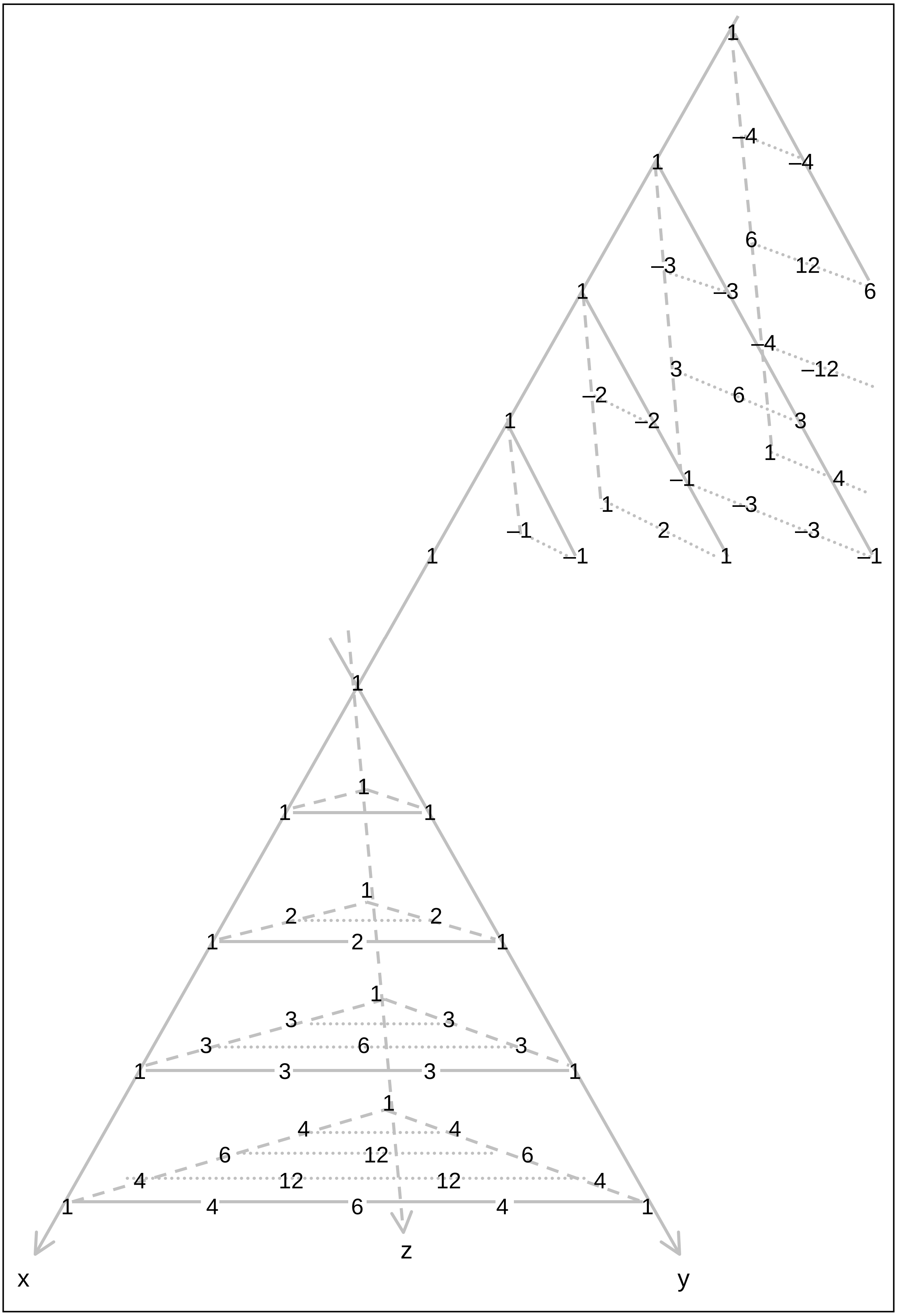

Fig. 3: Pascal space with positive and first negative Pascal pyramid.

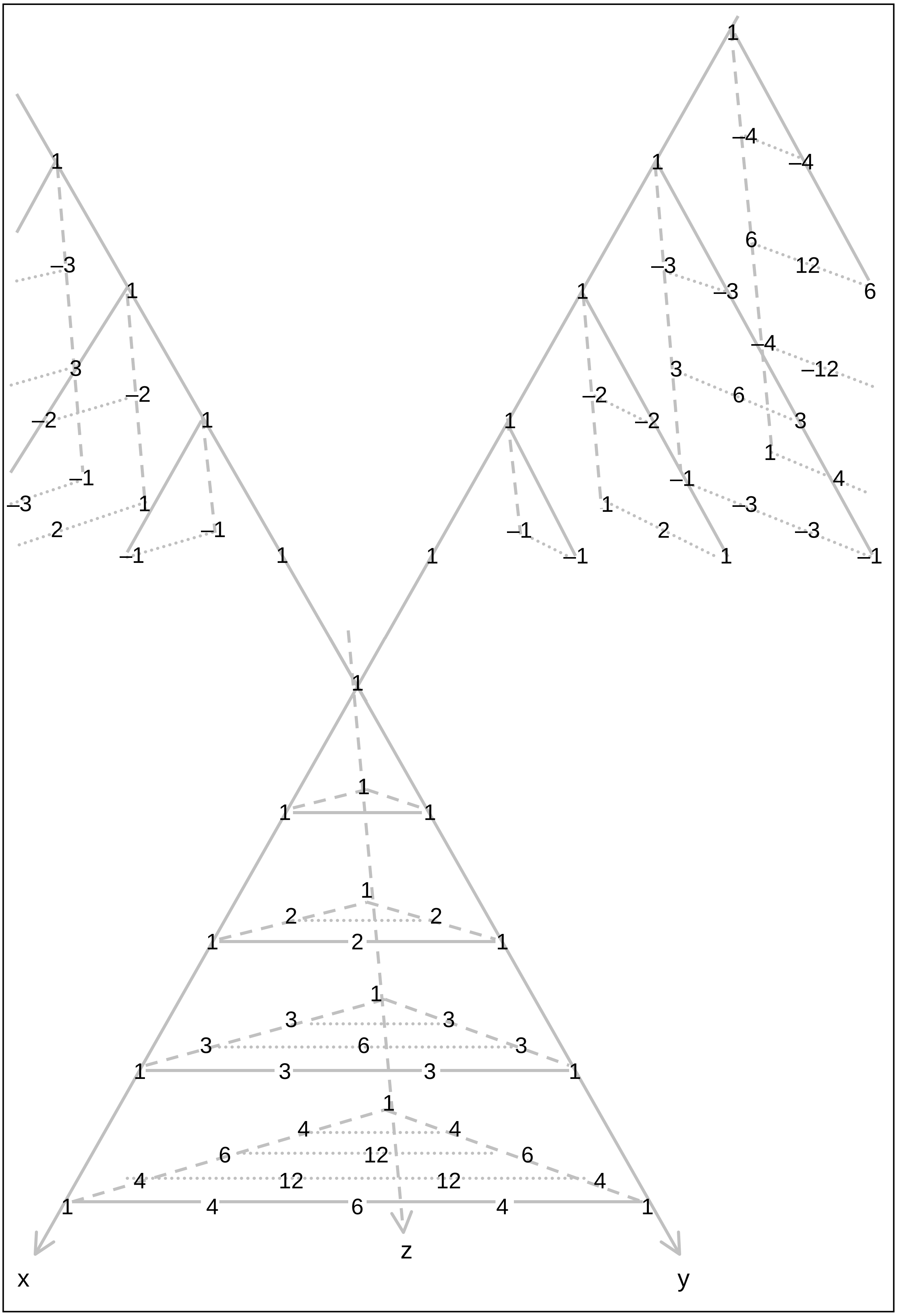

Fig. 4: Pascal space with positive, first negative, and second negative Pascal pyramid.

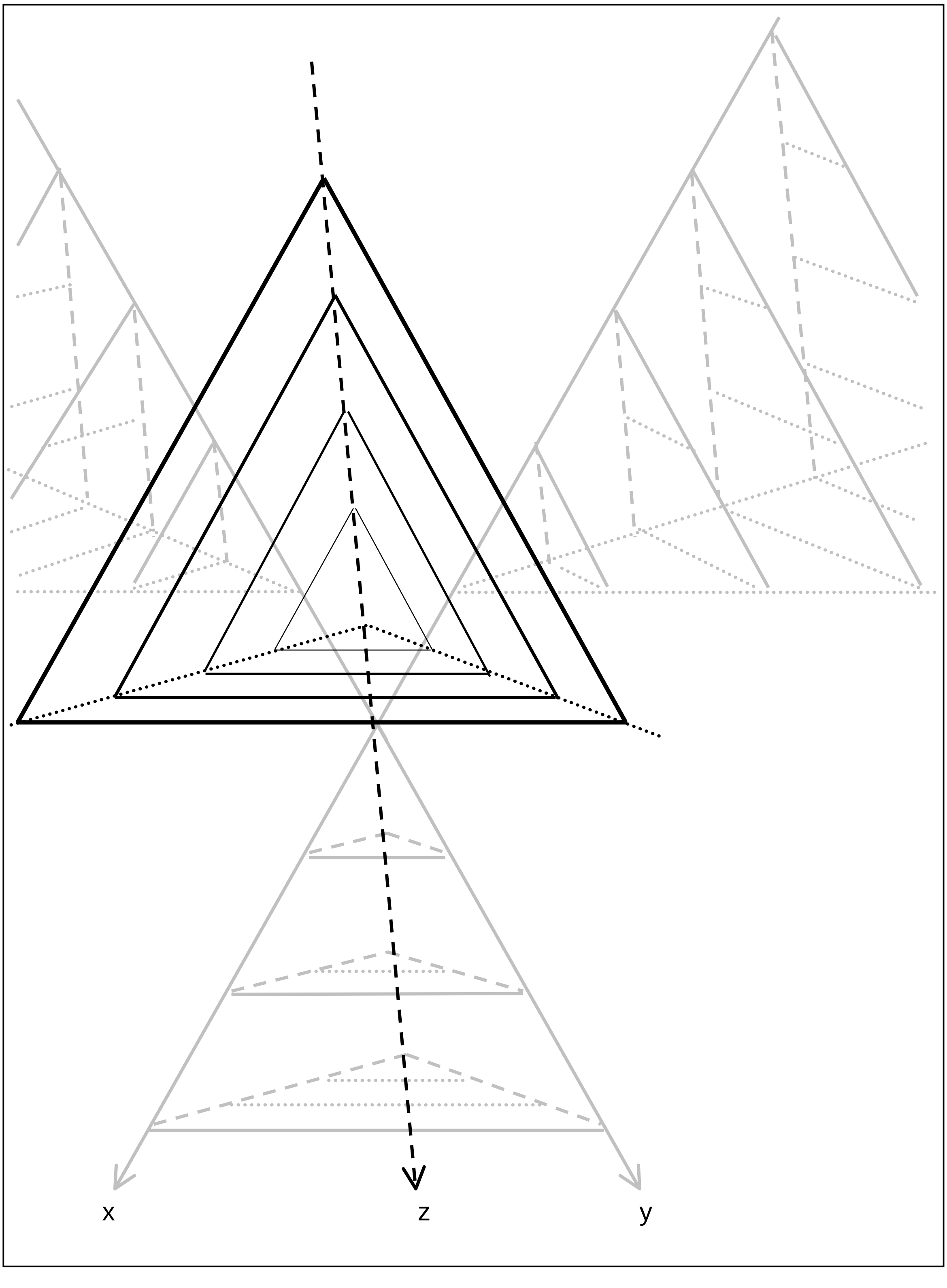

Fig. 5: The structure of Pascal space with all Pascal pyramids. The third negative Pascal pyramid is indicated by bold lines.
Fig. 12: The structure of Pauli Pascal space with all Pauli Pascal pyramids. The third negative Pauli Pascal pyramid is indicated by bold lines.

## 4. Basic Elements of Geometric Algebra

In the introduction some problems of contemporary mathematics like the unreadableness of mathematical texts were mentioned. Mathematical texts usually are dense and compact. But it comes worse: They seem to be not only difficult to understand because of their density regarding the content, but also because of being written in different mathematical languages. “In mathematics (…) one gets the impression sometimes that definitions are being created more rapidly than theorems” [4]. Therefore “a Babylon of mathematical tongues contributes to fragmentation of knowledge” [5]. Even “mathematicians in neighboring disciplines can hardly talk to one another” [4].

Oersted prize winner David Hestenes therefore proposes a unified language [4] which will help physicists and mathematicians to speak about our 3-dimensional world. This construction is based on the geometric and algebraic structure of 3-dimensional space (see figure 6).

| | | | |
|---|---|---|---|
| **1** | $\sigma_x, \sigma_y, \sigma_z$ | $\sigma_x\sigma_y, \sigma_y\sigma_z, \sigma_z\sigma_x$ | $\sigma_x\sigma_y\sigma_z$ |
| scalar | vectors | bivectors | trivector |
| 1 element without dimension | 3 orthonormal directions | 3 directed area element | 1 directed volume element |

Fig. 6: Basic elements of geometric algebra.

The basic elements form a linear space of 8 dimensions. Because the three vectors algebraically behave like Pauli matrices, these matrices get another interpretation: In geometric algebra Pauli matrices can be seen as vectors of our 3-dimensional world we lived in before relativity was invented by Einstein [3].

## 5. Anti-commutative Pauli Pascal Triangles

While the three Pascal triangles codify the coefficients of the binomial series, the Pauli Pascal triangles will codify the coefficients of an anticommuting binomial series. It seems that such Pauli Pascal triangles are of epistemological interest, too [7].

Let $a_x := a\sigma_x$ and $b_y := b\sigma_y$

For $n \geq 0$ there is

$$\begin{aligned}
(a_x + b_y)^0 &= 1 \\
(a_x + b_y)^1 &= 1a_x + 1b_y \\
(a_x + b_y)^2 &= 1a_x^2 + 0a_xb_y + 1b_y^2 \\
(a_x + b_y)^3 &= 1a_x^3 + 1a_x^2b_y + 1a_xb_y^2 + 1b_y^3 \\
(a_x + b_y)^4 &= 1a_x^4 + 0a_x^3b_y + 2a_x^2b_y^2 + 0a_xb_y^3 + 1b_y^4 \\
(a_x + b_y)^5 &= 1a_x^5 + 1a_x^4b_y + 2a_x^3b_y^2 + 2a_x^2b_y^3 + 1a_xb_y^4 + 1b_y^5 \\
(a_x + b_y)^6 &= 1a_x^6 + 0a_x^5b_y + 3a_x^4b_y^2 + 0a_x^3b_y^3 + 3a_x^2b_y^4 + 0a_xb_y^5 + 1b_y^6 \\
&\text{etc…}
\end{aligned} \tag{6}$$

because of the anti-commutative structure of the Pauli matrices $\sigma_x\sigma_y = -\sigma_y\sigma_x$. The coefficients can be arranged as anti-commutative Pauli Pascal triangle. It is interesting that this pattern shows an inherent threefold structure (see figure 7). Therefore we should pay attention when Taylor expanding functions of Pauli matrices. An inherent threefold structure[1] may automatically appear.

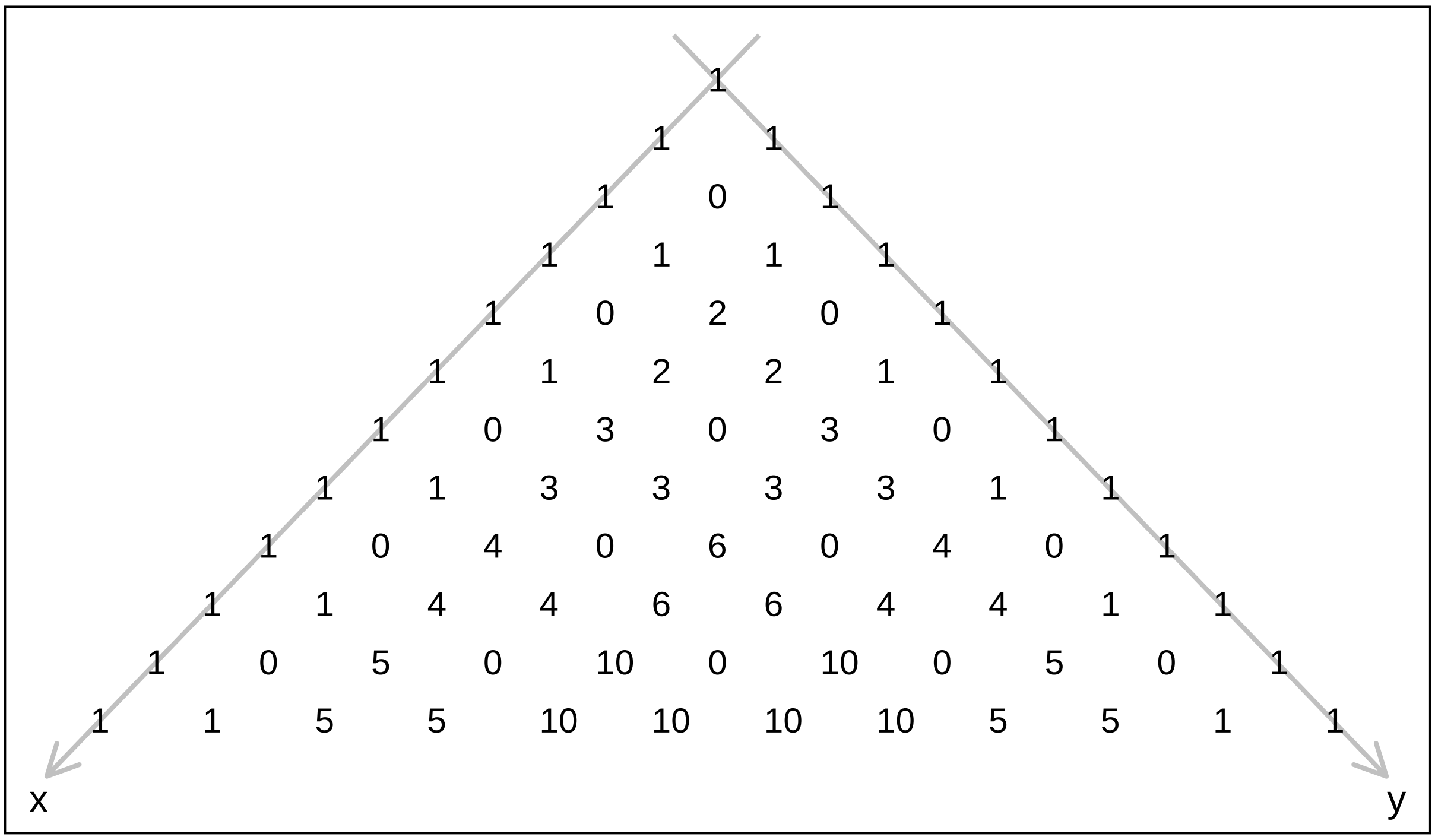

Fig. 7: The positive Pauli Pascal triangle.

Division through vectors or multivectors is not forbidden in geometric algebra. Therefore the two negative Pauli Pascal triangles can be determined as well.

For $|a| > |b|$ there is:

$$(a_x + b_y)^{-1} = \frac{a_x + b_y}{(a_x + b_y)^2} = \frac{a_x + b_y}{a_x^2 + b_y^2} = \frac{a_x + b_y}{a_x^2} \cdot \frac{1}{1 + b_y^2 / a_x^2} = \left(a_x^{-1} + a_x^{-2} b_y\right) \cdot \left(1 + \frac{b_y^2}{a_x^2}\right)^{-1} \quad (7a)$$

Because $(b_y^2/a_x^2)$ is a scalar the last term on the right hand side of (7a) can be expanded using the Taylor expansion $(1 + x) = 1 - x + x^2 - x^3 + x^4 - \ldots + \ldots$

$$(a_x + b_y)^{-1} = 1a_x^{-1} + 1a_x^{-2}b_y - 1a_x^{-3}b_y^2 - 1a_x^{-4}b_y^3 + 1a_x^{-5}b_y^4 + 1a_x^{-6}b_y^5 - \ldots \quad (7b)$$

$$(a_x + b_y)^{-2} = \frac{1}{(a_x + b_y)^2} = \frac{1}{a_x^2 + b_y^2} = a_x^{-2} \cdot \left(1 + \frac{b_y^2}{a_x^2}\right)^{-1}$$
$$= 1a_x^{-2} + 0a_x^{-3}b_y - 1a_x^{-4}b_y^2 + 0a_x^{-5}b_y^3 + 1a_x^{-6}b_y^4 + 0a_x^{-7}b_y^5 - \ldots \quad (8)$$

$$(a_x + b_y)^{-3} = (a_x + b_y)^{-2} (a_x + b_y)^{-1}$$
$$= 1a_x^{-3} + 1a_x^{-4}b_y - 2a_x^{-5}b_y^2 - 2a_x^{-6}b_y^3 + 3a_x^{-7}b_y^4 + 3a_x^{-8}b_y^5 - \ldots \quad (9)$$

$$(a_x + b_y)^{-4} = (a_x + b_y)^{-2} (a_x + b_y)^{-2}$$
$$= 1a_x^{-4} + 0a_x^{-5}b_y - 2a_x^{-6}b_y^2 + 0a_x^{-7}b_y^3 + 3a_x^{-8}b_y^4 + 0a_x^{-9}b_y^5 - \ldots \quad (10)$$

[1] Remember the quarks!

$$(a_x + b_y)^{-5} = (a_x + b_y)^{-4}\,(a_x + b_y)^{-1}$$
$$= 1a_x^{-5} + 1a_x^{-6}b_y - 3a_x^{-7}b_y^{2} - 3a_x^{-8}b_y^{3} + 6a_x^{-9}b_y^{4} + 6a_x^{-10}b_y^{5} - \ldots \qquad (11)$$

etc…

For $|a| < |b|$ there is:

$$
\begin{aligned}
(a_x + b_y)^{-1} &= 1a_x^{0}b_y^{-1} + 1a_x^{1}b_y^{-2} - 1a_x^{2}b_y^{-3} - 1a_x^{3}b_y^{-4} + 1a_x^{4}b_y^{-5} + 1a_x^{5}b_y^{-6} - \ldots \\
(a_x + b_y)^{-2} &= 1a_x^{0}b_y^{-2} + 0a_x^{1}b_y^{-3} - 1a_x^{2}b_y^{-4} + 0a_x^{3}b_y^{-5} + 1a_x^{4}b_y^{-6} + 0a_x^{5}b_y^{-7} - \ldots \\
(a_x + b_y)^{-3} &= 1a_x^{0}b_y^{-3} + 1a_x^{1}b_y^{-4} - 2a_x^{2}b_y^{-5} - 2a_x^{3}b_y^{-6} + 3a_x^{4}b_y^{-7} + 3a_x^{5}b_y^{-8} - \ldots \qquad (12) \\
(a_x + b_y)^{-4} &= 1a_x^{0}b_y^{-4} + 0a_x^{1}b_y^{-5} - 2a_x^{2}b_y^{-6} + 0a_x^{3}b_y^{-7} + 3a_x^{4}b_y^{-8} + 0a_x^{5}b_y^{-9} - \ldots \\
(a_x + b_y)^{-5} &= 1a_x^{0}b_y^{-5} + 1a_x^{1}b_y^{-6} - 3a_x^{2}b_y^{-7} - 3a_x^{3}b_y^{-8} + 6a_x^{4}b_y^{-9} + 6a_x^{5}b_y^{-10} - \ldots
\end{aligned}
$$

etc…

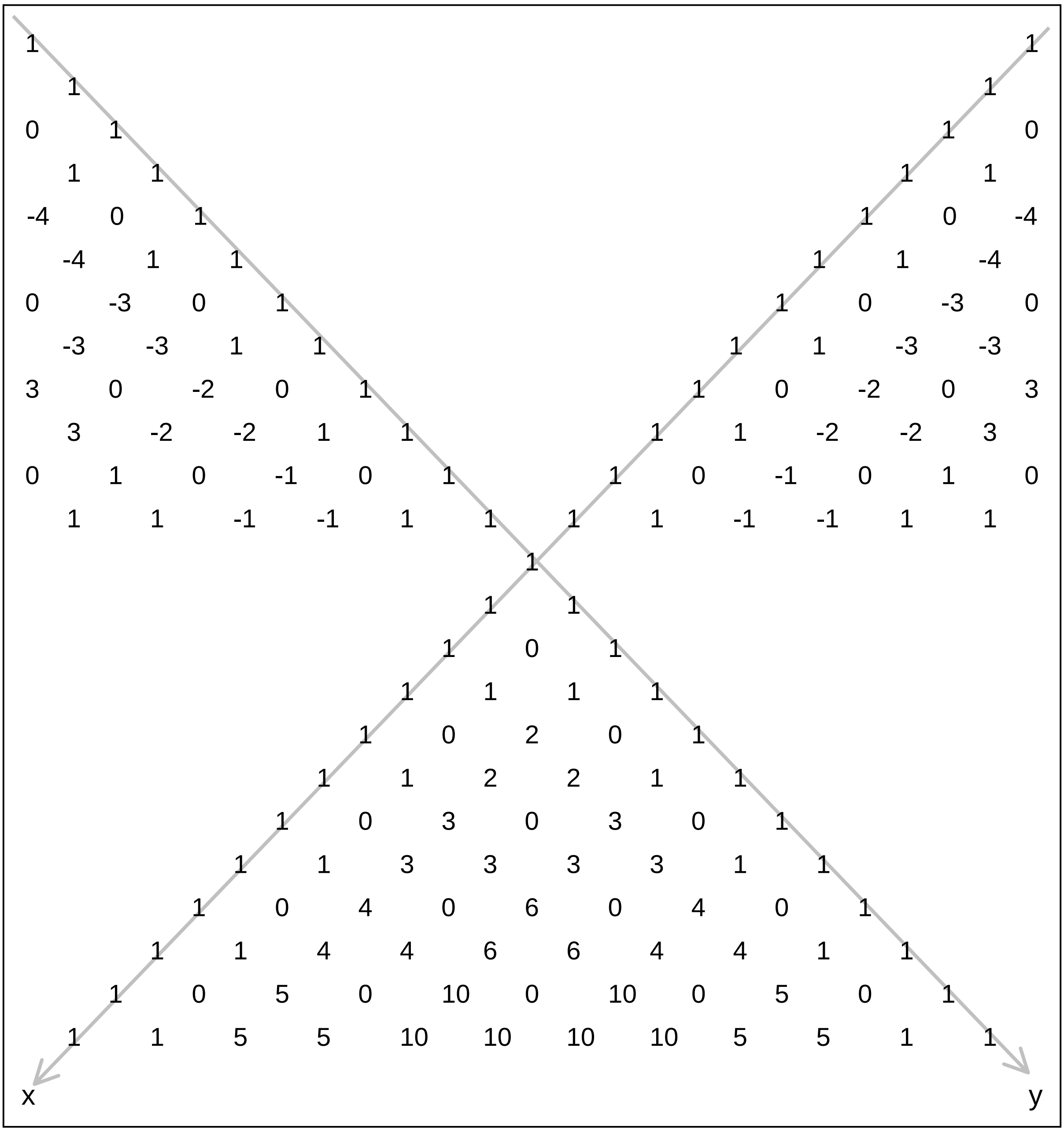

Fig. 8: The positive Pauli Pascal triangle and the two negative Pauli Pascal triangles.

As figure 8 shows, all three Pauli Pascal triangles have an inherent threefold structure. When Taylor expanding functions of Pauli matrices this inherent threefold structure[2] may automatically appear three times.

## 6. Pauli Pascal Pyramids and Pauli Pascal Hyperpyramids

Let $a_x := a\sigma_x$, $b_y := b\sigma_y$, and $c_z := c\sigma_z$

For $n \geq 0$ there is

$$
\begin{aligned}
(a_x + b_y + c_z)^0 &= 1 \\
(a_x + b_y + c_z)^1 &= (1a_x + 1b_y) + 1c_z \\
(a_x + b_y + c_z)^2 &= (1a_x^2 + 0a_xb_y + 1b_y^2) + (0a_xc_z + 0b_yc_z) + 1c_z^2 \\
(a_x + b_y + c_z)^3 &= (1a_x^3 + 1a_x^2b_y + 1a_xb_y^2 + 1b_y^3) + (1a_x^2c_z + 0a_xb_yc_z + 1b_y^2c_z) \\
&\quad + (1\,a_xc_z^2 + 1b_yc_z^2) + 1c_z^3 \\
(a_x + b_y + c_z)^4 &= (1a_x^4 + 0a_x^3b_y + 2a_x^2b_y^2 + 0a_xb_y^3 + 1b_y^4) \\
&\quad + (0a_x^3c_z + 0a_x^2b_yc_z + 0a_xb_y^2c_z + 0b_y^3c_z) \\
&\quad + (2a_x^2c_z^2 + 0a_xb_yc_z^2 + 2b_y^2c_z^2) + (0a_xc_z^3 + 0b_yc_z^3) + 1c_z^4 \\
(a_x + b_y + c_z)^5 &= (1a_x^5 + 1a_x^4b_y + 2a_x^3b_y^2 + 2a_x^2b_y^3 + 1a_xb_y^4 + 1b_y^5) \qquad (13) \\
&\quad + (1a_x^4c_z + 0a_x^3b_yc_z + 2a_x^2b_y^2c_z + 0a_xb_y^3c_z + 1b_y^4c_z) \\
&\quad + (2a_x^3c_z^2 + 2a_x^2b_yc_z^2 + 2a_xb_y^2c_z^2 + 2b_y^3c_z^2) \\
&\quad + (2a_x^2c_z^3 + 0a_xb_yc_z^3 + 2b_y^2c_z^3) + (1a_xc_z^4 + 1b_yc_z^4) + 1c_z^5 \\
(a_x + b_y + c_z)^6 &= (1a_x^6 + 0a_x^5b_y + 3a_x^4b_y^2 + 0a_x^3b_y^3 + 3a_x^2b_y^4 + 0a_xb_y^5 + 1b_y^6) \\
&\quad + (0a_x^5c_z + 0a_x^4b_yc_z + 0a_x^3b_y^2c_z + 0a_x^2b_y^3c_z + 0a_xb_y^4c_z + 0b_y^5c_z) \\
&\quad + (3a_x^4c_z^2 + 0a_x^3b_yc_z^2 + 6a_x^2b_y^2c_z^2 + 0a_xb_y^3c_z^2 + 3b_y^4c_z^2) \\
&\quad + (0a_x^3c_z^3 + 0a_x^2b_yc_z^3 + 0a_xb_y^2c_z^3 + 0b_y^3c_z^3) \\
&\quad + (3a_x^2c_z^4 + 0a_xb_yc_z^4 + 3b_y^2c_z^4) + (0a_xc_z^5 + 0b_yc_z^5) + 1c_z^6
\end{aligned}
$$

etc…

Again the coefficients can be arranged as Pauli Pascal pyramid (see figure 9).

For $n < 0$ and $|\,a\,| > |\,b\,|, |\,c\,|$ by using formulae (5) there is:

$$
\begin{aligned}
&(a_x + b_y + c_z)^{-1} \\
&\quad = (a_x + b_y + c_z)\,(a_x + b_y + c_z)^{-2} \\
&\quad = (a_x + b_y + c_z)\,(a_x^2 + b_y^2 + c_z^2)^{-1} \\
&\quad = (a_x + b_y + c_z)\,(1a_x^{-2} - 1a_x^{-4}b_y^2 - 1a_x^{-4}c_z^2 + 1a_x^{-6}b_y^4 + 2a_x^{-6}b_y^2c_z^2 + 1a_x^{-6}c_z^4 \\
&\qquad\qquad - 1a_x^{-8}b_y^6 - 3a_x^{-8}b_y^4c_z^2 - 3a_x^{-8}b_y^2c_z^4 - 1a_x^{-8}c_z^6 + \ldots - \ldots) \\
&\quad = 1a_x^{-1} + 1a_x^{-2}b_y + 1a_x^{-2}c_z - 1a_x^{-3}b_y^2 - 1a_x^{-3}c_z^2 \qquad (14) \\
&\qquad - 1a_x^{-4}b_y^3 - 1a_x^{-4}b_y^2c_z - 1a_x^{-4}b_yc_z^3 - 1a_x^{-4}c_z^3 \\
&\qquad + 1a_x^{-5}b_y^4 + 2a_x^{-5}b_y^2c_z^2 + 1a_x^{-5}c_z^4 \\
&\qquad + 1a_x^{-6}b_y^5 + 1a_x^{-6}b_y^4c_z + 2a_x^{-6}b_y^3c_z^2 + 2a_x^{-6}b_y^2c_z^3 + 1a_x^{-6}b_yc_z^4 + 1a_x^{-6}c_z^5 \\
&\qquad - 1a_x^{-7}b_y^6 - 3a_x^{-7}b_y^4c_z^2 - 3a_x^{-7}b_y^2c_z^4 - 1a_x^{-7}c_z^6 + \ldots - \ldots
\end{aligned}
$$

[2] Remember the quarks! And don't forget rishons or preons…

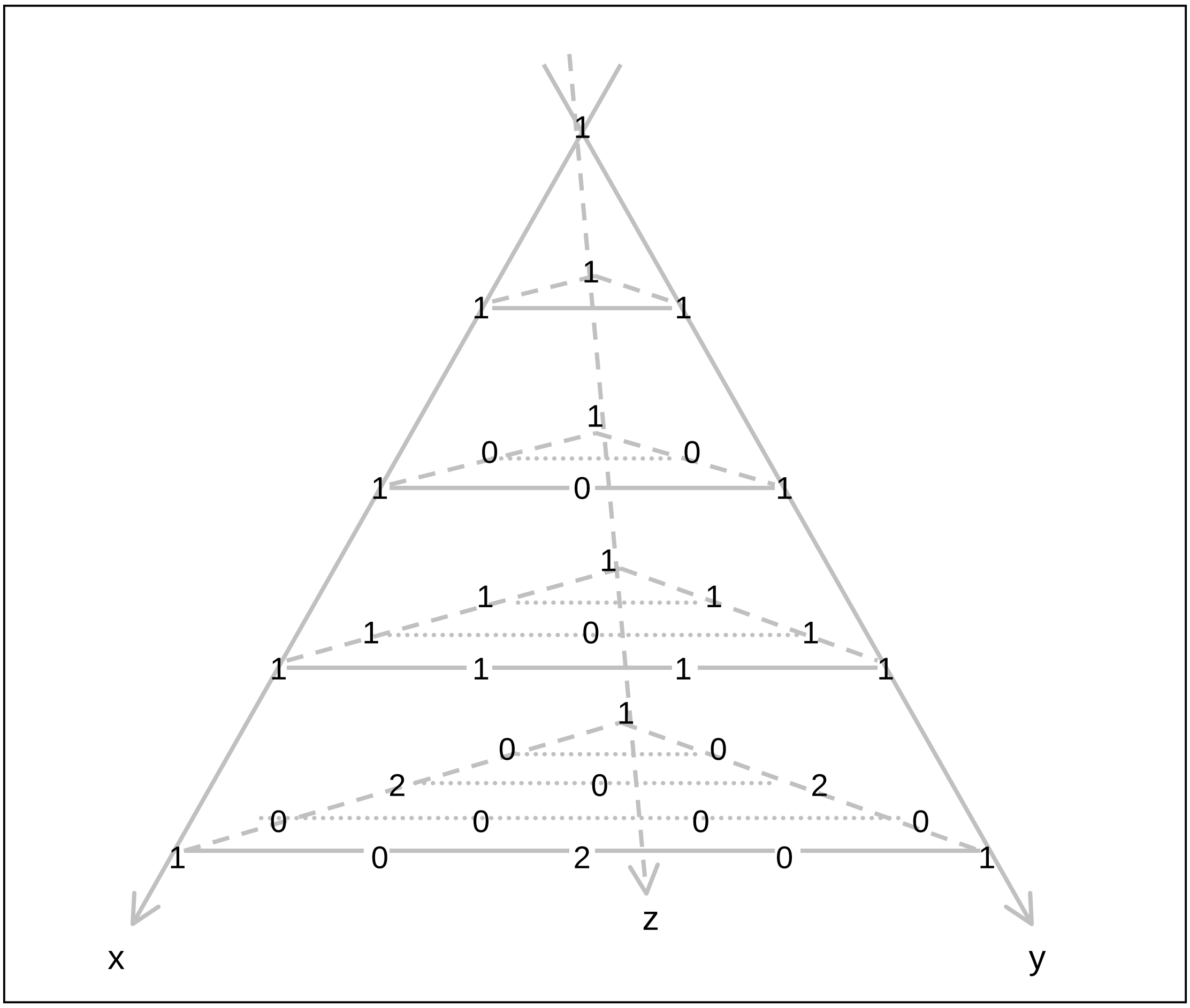

Fig. 9: Pauli Pascal space with the positive Pauli Pascal pyramid.

$$
\begin{aligned}
&(a_x + b_y + c_z)^{-2} \\
&\quad = (a_x^2 + b_y^2 + c_z^2)^{-1} \\
&\quad = 1a_x^{-2} - 1a_x^{-4}b_y^2 - 1a_x^{-4}c_z^2 + 1a_x^{-6}b_y^4 + 2a_x^{-6}b_y^2c_z^2 + 1a_x^{-6}c_z^4 \\
&\qquad - 1a_x^{-8}b_y^6 - 3a_x^{-8}b_y^4c_z^2 - 3a_x^{-8}b_y^2c_z^4 - 1a_x^{-8}c_z^6 + \ldots - \ldots
\end{aligned} \tag{15}
$$

$$
\begin{aligned}
&(a_x + b_y + c_z)^{-3} \\
&\quad = (a_x + b_y + c_z)\,(a_x + b_y + c_z)^{-4} \\
&\quad = (a_x + b_y + c_z)\,(a_x^2 + b_y^2 + c_z^2)^{-2} \\
&\quad = (a_x + b_y + c_z)\,(1a_x^{-4} - 2a_x^{-6}b_y^2 - 2a_x^{-6}c_z^2 + 3a_x^{-8}b_y^4 + 6a_x^{-8}b_y^2c_z^2 + 3a_x^{-8}c_z^4 \\
&\qquad - 4a_x^{-10}b_y^6 - 12a_x^{-10}b_y^4c_z^2 - 12a_x^{-10}b_y^2c_z^4 - 4a_x^{-10}c_z^6 + \ldots - \ldots) \\
&\quad = 1a_x^{-3} + 1a_x^{-4}b_y + 1a_x^{-4}c_z - 2a_x^{-5}b_y^2 - 2a_x^{-5}c_z^2 \\
&\qquad - 2a_x^{-6}b_y^3 - 2a_x^{-6}b_y^2c_z - 2a_x^{-6}b_yc_z^2 - 2a_x^{-6}c_z^3 + 3a_x^{-7}b_y^4 + 6a_x^{-7}b_y^2c_z^2 + 3a_x^{-7}c_z^4 \\
&\qquad + 3a_x^{-8}b_y^5 + 3a_x^{-8}b_y^4c_z + 6a_x^{-8}b_y^3c_z^2 + 6a_x^{-8}b_y^2c_z^3 + 3a_x^{-8}b_yc_z^4 + 3a_x^{-8}c_z^5 \\
&\qquad - 4a_x^{-9}b_y^6 - 12a_x^{-9}b_y^4c_z^2 - 12a_x^{-9}b_y^2c_z^4 - 4a_x^{-9}c_z^6 - \ldots + \ldots
\end{aligned} \tag{16}
$$

$$
\begin{aligned}
&(a_x + b_y + c_z)^{-4} \\
&\quad = (a_x^2 + b_y^2 + c_z^2)^{-2} \\
&\quad = 1a_x^{-4} - 2a_x^{-6}b_y^2 - 2a_x^{-6}c_z^2 + 3a_x^{-8}b_y^4 + 6a_x^{-8}b_y^2c_z^2 + 3a_x^{-8}c_z^4 \\
&\qquad - 4a_x^{-10}b_y^6 - 12a_x^{-10}b_y^4c_z^2 - 12a_x^{-10}b_y^2c_z^4 - 4a_x^{-10}c_z^6 + \ldots - \ldots
\end{aligned} \tag{17}
$$

$$
\begin{aligned}
&(a_x + b_y + c_z)^{-5} \\
&\quad = (a_x + b_y + c_z)\,(a_x + b_y + c_z)^{-6} \\
&\quad = (a_x + b_y + c_z)\,(a_x^2 + b_y^2 + c_z^2)^{-3} \\
&\quad = (a_x + b_y + c_z)\,(1a_x^{-6} - 3a_x^{-8}b_y^2 - 3a_x^{-8}c_z^2 + 6a_x^{-10}b_y^4 + 12a_x^{-10}b_y^2c_z^2 + 6a_x^{-10}c_z^4 \\
&\qquad\qquad - 10a_x^{-12}b_y^6 - 30a_x^{-12}b_y^4c_z^2 - 30a_x^{-12}b_y^2c_z^4 - 10a_x^{-12}c_z^6 + \ldots - \ldots) \\
&\quad = 1a_x^{-5} + 1a_x^{-6}b_y + 1a_x^{-6}c_z - 3a_x^{-7}b_y^2 - 3a_x^{-7}c_z^2 \qquad (18) \\
&\qquad - 3a_x^{-8}b_y^3 - 3a_x^{-8}b_y^2c_z - 3a_x^{-8}b_yc_z^2 - 3a_x^{-8}c_z^3 + 6a_x^{-9}b_y^4 + 12a_x^{-9}b_y^2c_z^2 + 6a_x^{-9}c_z^4 \\
&\qquad + 6a_x^{-10}b_y^5 + 6a_x^{-10}b_y^4c_z + 12a_x^{-10}b_y^3c_z^2 + 12a_x^{-10}b_y^2c_z^3 + 6a_x^{-10}b_yc_z^4 + 6a_x^{-10}c_z^5 \\
&\qquad - 10a_x^{-11}b_y^6 - 30a_x^{-11}b_y^4c_z^2 - 30a_x^{-11}b_y^2c_z^4 - 10a_x^{-11}c_z^6 - \ldots + \ldots
\end{aligned}
$$

$$
\begin{aligned}
&(a_x + b_y + c_z)^{-6} \\
&\quad = (a_x^2 + b_y^2 + c_z^2)^{-3} \\
&\quad = 1a_x^{-6} - 3a_x^{-8}b_y^2 - 3a_x^{-8}c_z^2 + 6a_x^{-10}b_y^4 + 12a_x^{-10}b_y^2c_z^2 + 6a_x^{-10}c_z^4 \qquad (19) \\
&\qquad - 10a_x^{-12}b_y^6 - 30a_x^{-12}b_y^4c_z^2 - 30a_x^{-12}b_y^2c_z^4 - 10a_x^{-12}c_z^6 + \ldots - \ldots
\end{aligned}
$$

These coefficients again can be arranged as first negative Pauli Pascal pyramid (see figure 10). The other two negative Pauli Pascal pyramids can be constructed in a similar way. If the coefficients are arranged in the usual way, the second negative Pauli Pascal pyramid (see figure 11) and the third negative Pauli Pascal pyramid (see figure 12) will emerge. While the surface of the positive Pauli Pascal pyramid is composed of the positive Pauli Pascal triangle, the surfaces of the negative Pauli Pascal pyramids are composed of negative Pascal triangles.

In four-dimensional Pauli Pascal hyperspace there exist one positive four-dimensional Pauli Pascal hyperpyramid and four negative four-dimensional Pauli Pascal hyperpyramids. The positive Pauli Pascal hyperpyramid can be constructed using

$$a_x := a\sigma_x = a\sigma_{x1}, \quad b_y := b\sigma_y = b\sigma_{x2}, \quad c_z := c\sigma_z = c\sigma_{x3}, \quad \text{and} \quad d_w := d\sigma_w = d\sigma_{x4},$$

with four independent Pauli matrices $\sigma_{xi}\sigma_{xj} = -\sigma_{xj}\sigma_{xi}$ for $i \neq j$.

The coefficients of the Taylor expansion of

$$(a_x + b_y + c_z + d_w)^n \quad \text{with } n \geq 0$$

or alternatively multinomial coefficients [8] can be arranged in four-dimensional Pauli Pascal space. The four three-dimensional hypersurfaces of this positive four-dimensional Pauli Pascal hyperpyramid are the positive Pauli Pascal pyramid of figure 9.

Negative Pauli Pascal hyperpyramids can be constructed with the coefficients of the Taylor expansion of $(a + b + c + d)^n$ with $n < 0$. The four three-dimensional hypersurfaces of these negative four-dimensional Pauli Pascal hyperpyramids are the negative Pauli Pascal pyramids of figure 11.

In a similar way $(k + 1)$ Pauli Pascal hyperpyramids of dimension k can be constructed with the coefficients of the Taylor expansion of $(a_1 + a_2 + a_3 + \ldots + a_k)^n$.

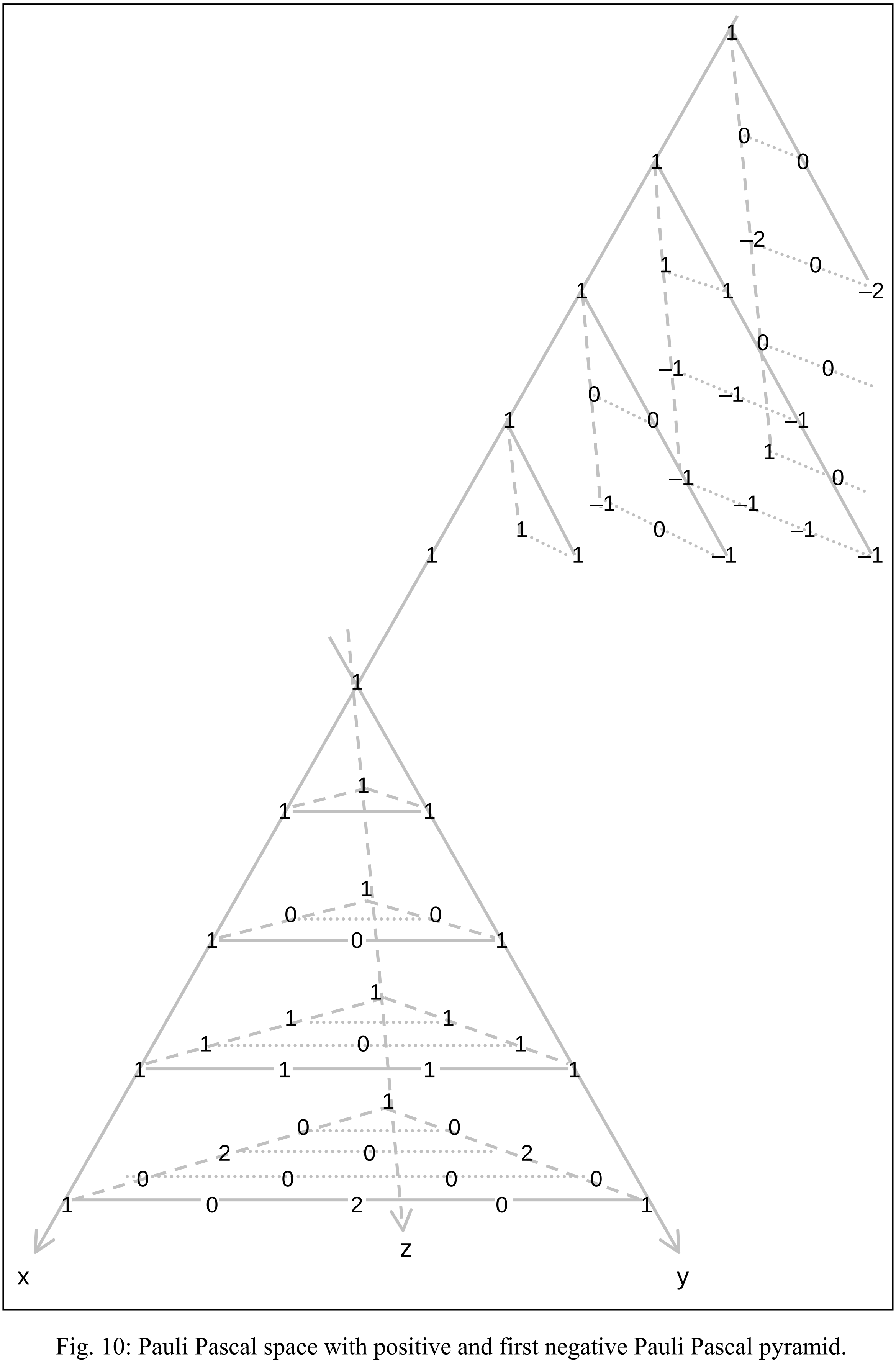

Fig. 10: Pauli Pascal space with positive and first negative Pauli Pascal pyramid.

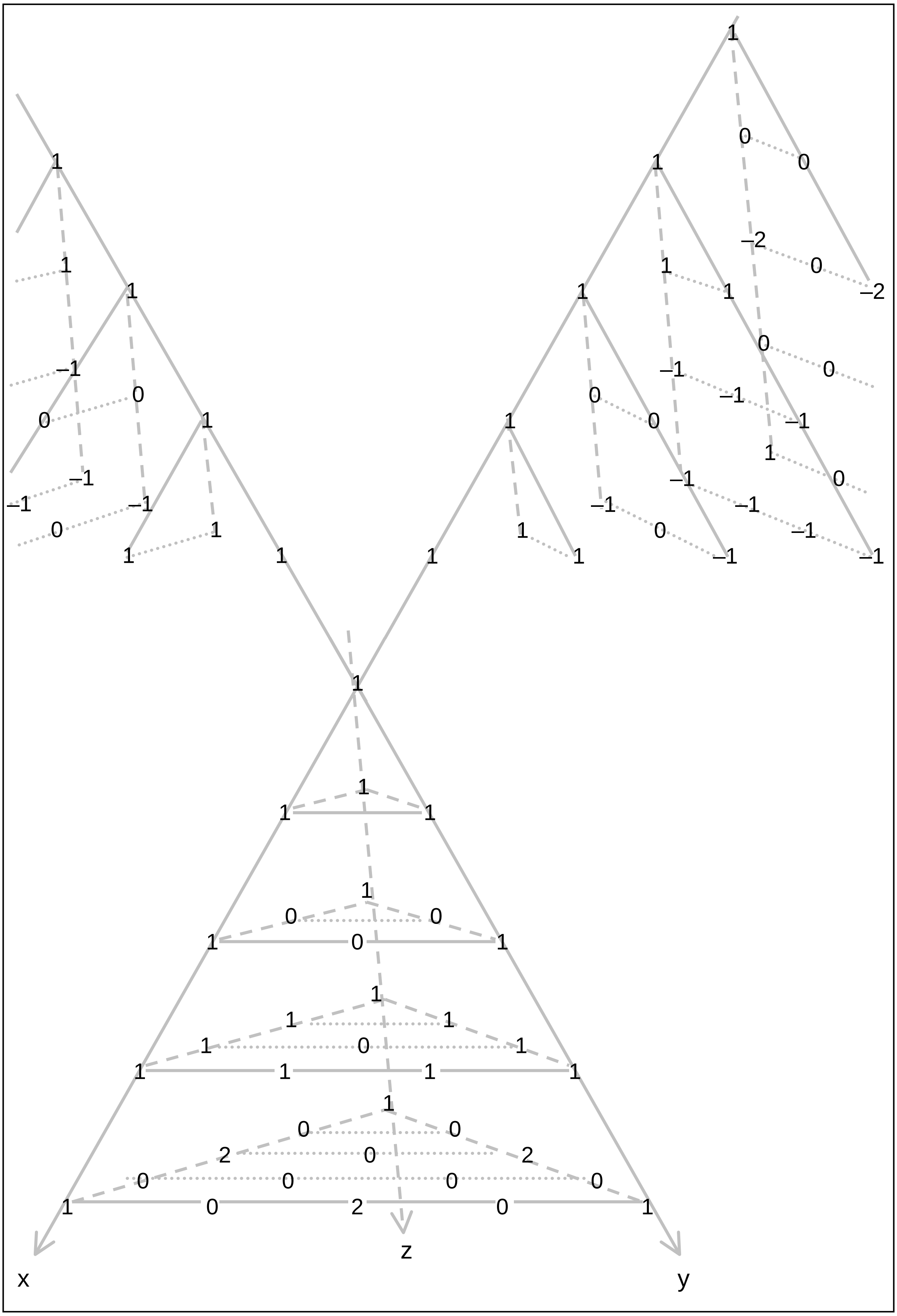

Fig. 11: Pauli Pascal space with positive, first negative, and the second negative Pauli Pascal pyramid.

## 7. Positive and Negative Fibonacci Numbers

If every pair of rabbits produces one pair of rabbits of the next generation, they can decode their Fibonaccian breeding habits with the help of the Pascal triangle (see figure 13).

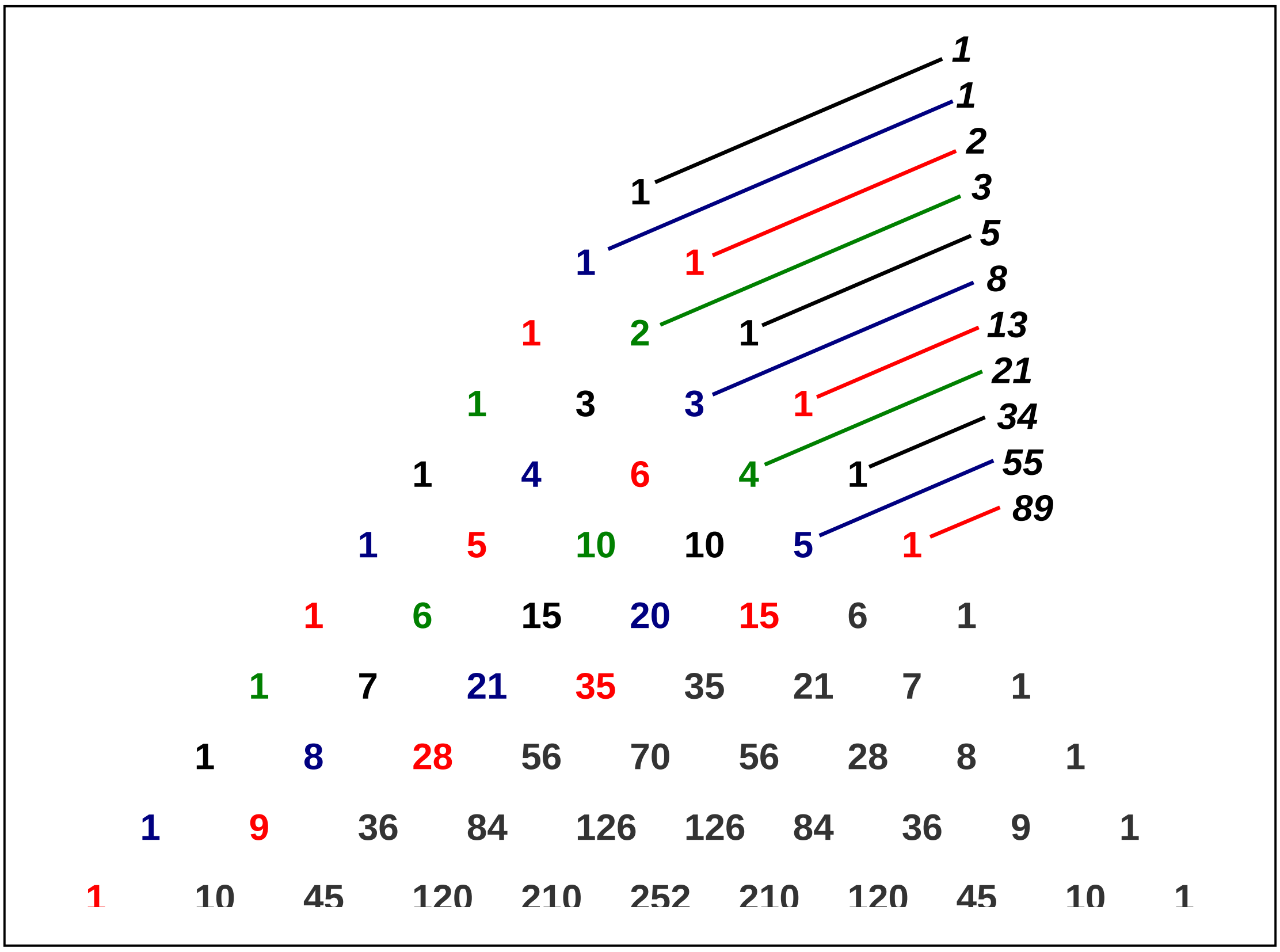

Fig. 13: Positive half of Fibonacci numbers.

They satisfy the recurrence relation

$$F_n = F_{n-1} + F_{n-2} \tag{20}$$

It is interesting to note, that Binet's formula

$$F_n = \frac{1}{\sqrt{5}}\left[\left(\frac{1+\sqrt{5}}{2}\right)^n - \left(\frac{1-\sqrt{5}}{2}\right)^n\right] \tag{21}$$

is also valid for $n \le 0$:

$$F_{-n} = \frac{1}{\sqrt{5}}\left[\left(\frac{1+\sqrt{5}}{2}\right)^{-n} - \left(\frac{1-\sqrt{5}}{2}\right)^{-n}\right] = \frac{1}{\sqrt{5}}\left[\left(\frac{2}{1+\sqrt{5}}\right)^n - \left(\frac{2}{1-\sqrt{5}}\right)^n\right]$$

$$= \frac{1}{\sqrt{5}}\left[\left(\frac{2\left(1-\sqrt{5}\right)}{1-5}\right)^n - \left(\frac{2\left(1+\sqrt{5}\right)}{1-5}\right)^n\right] = \frac{1}{\sqrt{5}}\left[\left(-\frac{1-\sqrt{5}}{2}\right)^n - \left(-\frac{1+\sqrt{5}}{2}\right)^n\right]$$

$$= \frac{(-1)^{n-1}}{\sqrt{5}} \left[ \left( \frac{1+\sqrt{5}}{2} \right)^{n} - \left( \frac{1-\sqrt{5}}{2} \right)^{n} \right] = (-1)^{n-1} \cdot F_n \tag{22}$$

Of course there is a negative extension of Fibonacci numbers (see figure 14).

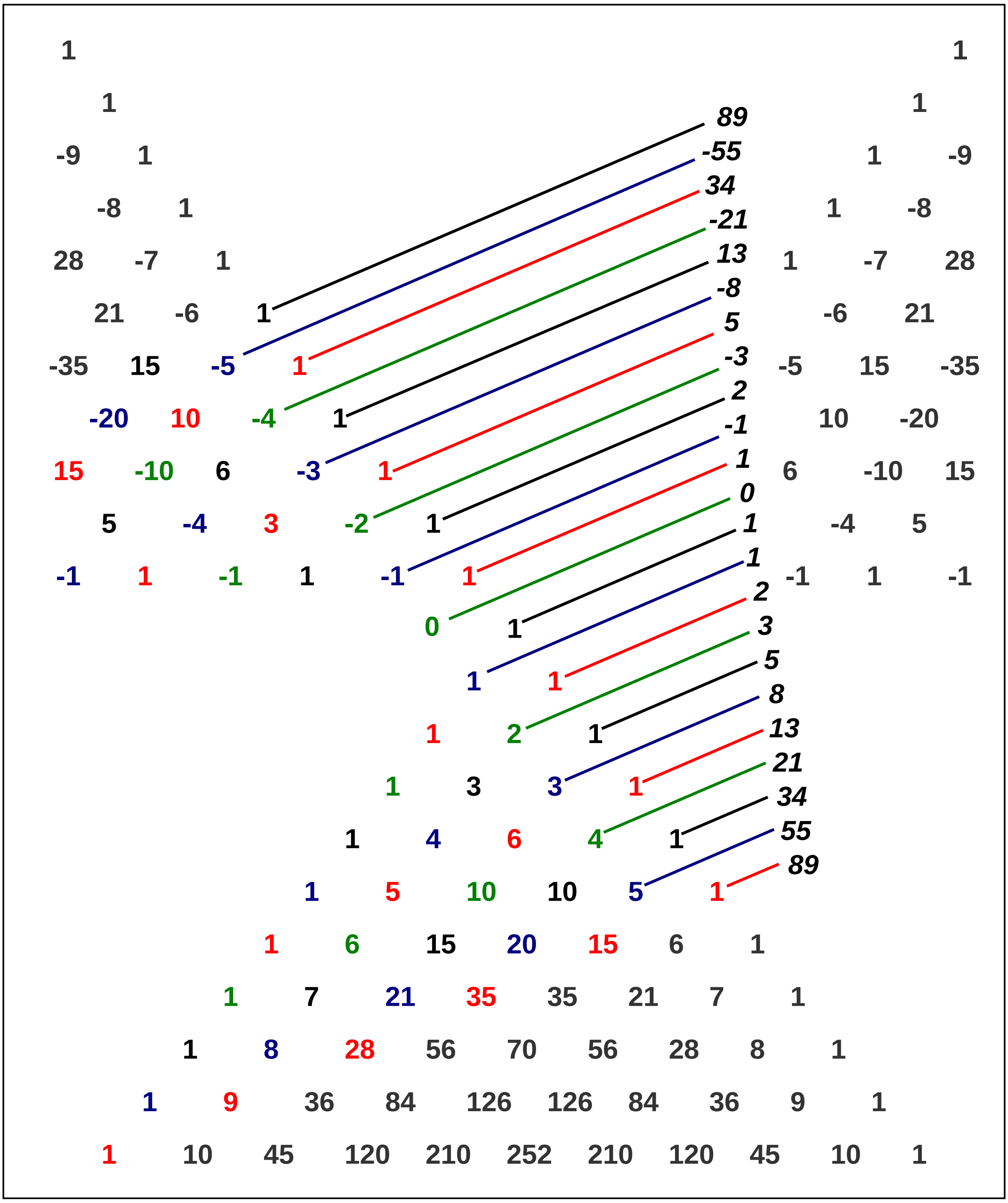

Fig. 14: Negative and positive Fibonacci numbers.

## 8. Pauli Fibonacci Numbers

If every pair of Pauli rabbits produces one pair of Pauli rabbits of the next generation, they can decode their Fibonaccian breeding habits with the help of the Pauli Pascal triangle (see figure 15).

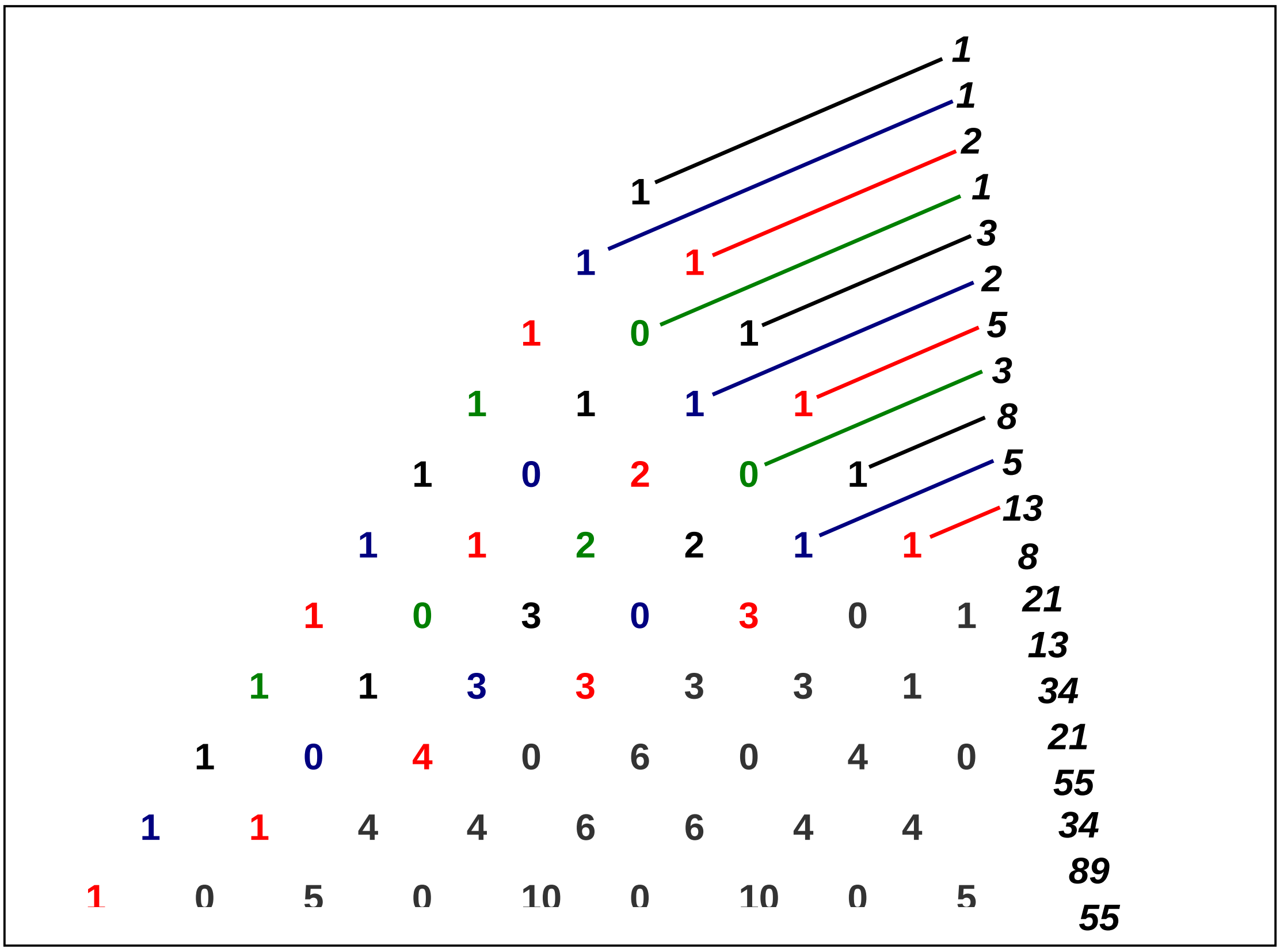

Fig. 15: Positive Pauli Fibonacci numbers.

There are two possibilities to read and to understand this result. First possibility: there are two totally distinct populations of rabbits which have nothing to do with each other and live in total isolation (see figure 16).

| | *period 1* | *period 2* | *period 3* | *period 4* | *period 5* | *period 6* | *period 7* | *period 8* | *period 9* | *period 10* |
|---|---|---|---|---|---|---|---|---|---|---|
| *first population of rabbits* | 1 | | 2 | | 3 | | 5 | | 8 | |
| *second population of rabbits* | | 1 | | 1 | | 2 | | 3 | | 5 |

Fig. 16: Two populations of rabbits living in total isolation.

Second possibility: A young pair of Pauli rabbits gives birth as usual to a pair of Pauli rabbits. But every second generation the parents then develop into killer rabbits and kill a pair of rabbits of the next generation (see figure 17).

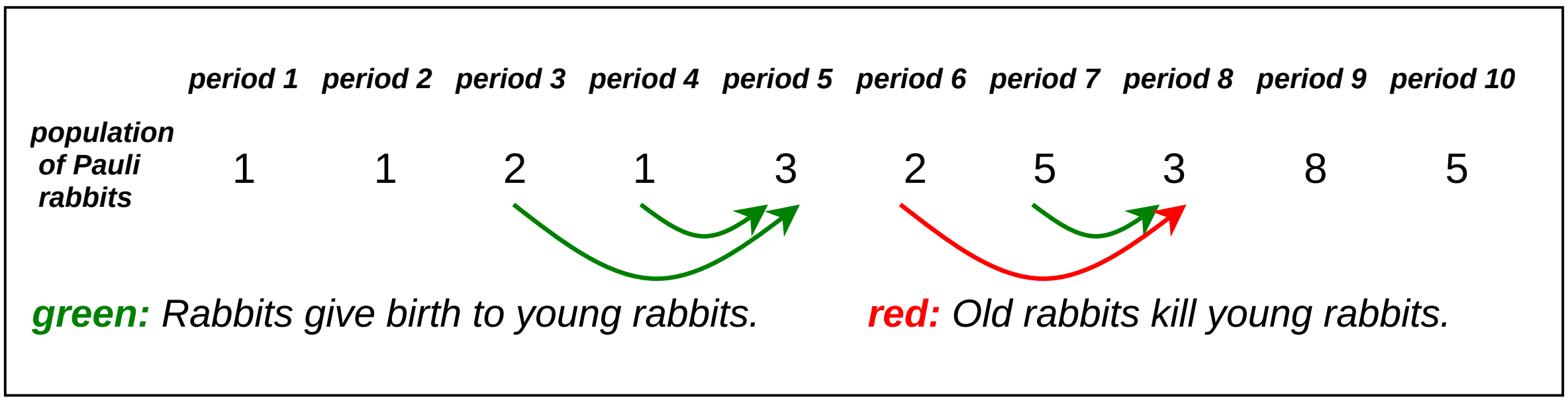

Fig. 17: Pauli rabbits are killer rabbits every second period.

If Pauli matrices are used in theoretical physics to model physical phenomena mathematically these two interpretations of Pauli matrices will play an important role. It is then possible to see a sort of period doubling with a twofold increase of the period length (first possibility). Or it is possible to see sort of interference effects with constructive and destructive interference (second possibility). Some problems will arise when both interpretations are mixed.

These two interpretations lead to two strategies, to modify Binet's formula for Pauli Fibonacci numbers $PF_n$. One possibility is, according to figure 17, to consider a modification of Binet's formula as a whole. Pauli Fibonacci numbers then can, for example, be generated by

$$PF_n = \frac{1}{\sqrt{5}}\left[\left(\frac{1+\sqrt{5}}{2}\right)^{\frac{2n+3+3\cos(\pi(n-1))}{4}} - \left(\frac{1-\sqrt{5}}{2}\right)^{\frac{2n+3+3\cos(\pi(n-1))}{4}}\right] \qquad (23)$$

with $n \in \mathbf{Z}$.

Another possibility is, according to figure 16, to split the Pauli Fibonacci numbers into two parts as in (24) and to modify every part separately.

$$\begin{array}{cccccccccccc} PF_{-1} = F_1 & & PF_1 = F_2 & & PF_3 = F_3 & & PF_5 = F_4 & & PF_7 = F_5 & & PF_9 = F_6 & \\ & PF_0 = F_0 & & PF_2 = F_1 & & PF_4 = F_2 & & PF_6 = F_3 & & PF_8 = F_4 & & \end{array} \qquad (24)$$

The second row of (24) can simply be generated by

$$PF_n^{(2)} = \frac{1}{2\sqrt{5}}\left[\sqrt{\frac{1+\sqrt{5}}{2}}^{\,n} + \left(-\sqrt{\frac{1+\sqrt{5}}{2}}\right)^{n} - \sqrt{\frac{1-\sqrt{5}}{2}}^{\,n} - \left(-\sqrt{\frac{1-\sqrt{5}}{2}}\right)^{n}\right] \qquad (25)$$

leaving space for the first row of (24) at odd n. This first row can be generated by

$$PF_n^{(1)} = \frac{1}{2\sqrt{5}}\left[\sqrt{\frac{1+\sqrt{5}}{2}}^{\,n+3} + \left(-\sqrt{\frac{1+\sqrt{5}}{2}}\right)^{n+3} - \sqrt{\frac{1-\sqrt{5}}{2}}^{\,n+3} - \left(-\sqrt{\frac{1-\sqrt{5}}{2}}\right)^{n+3}\right] \qquad (26)$$

The sum of (25) and (26) gives a Binet-like formula for Pauli Fibonacci numbers

$$PF_n = \frac{1}{2\sqrt{5}}\left[\left(1+\sqrt{2+\sqrt{5}}\right)\sqrt{\frac{1+\sqrt{5}}{2}}^{\,n} + \left(1-\sqrt{2+\sqrt{5}}\right)\left(-\sqrt{\frac{1+\sqrt{5}}{2}}\right)^{n} - \left(1+\sqrt{2-\sqrt{5}}\right)\sqrt{\frac{1-\sqrt{5}}{2}}^{\,n} - \left(1-\sqrt{2-\sqrt{5}}\right)\left(-\sqrt{\frac{1-\sqrt{5}}{2}}\right)^{n}\right] \tag{27}$$

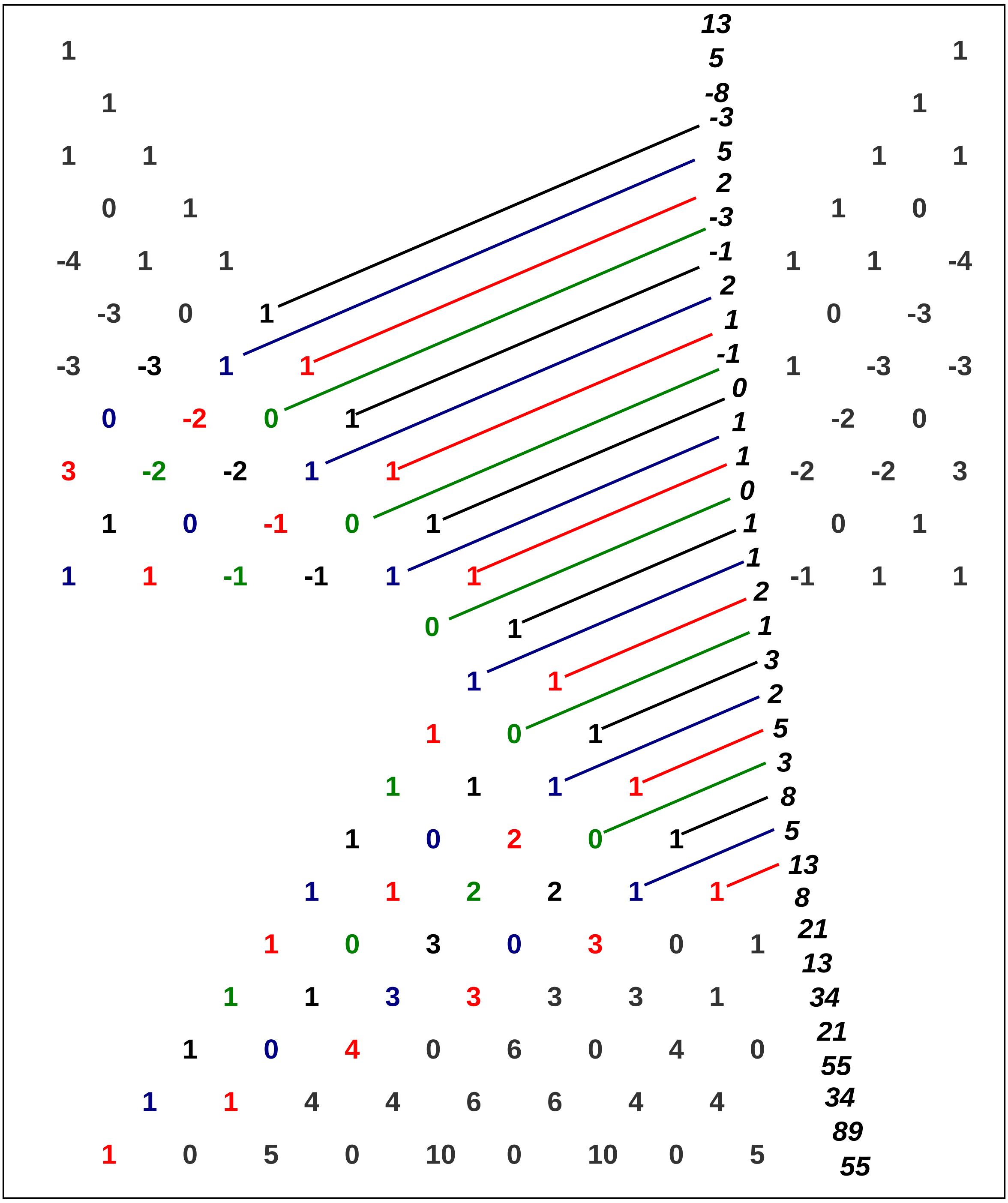

Fig. 18: Negative and positive Pauli Fibonacci numbers.

According to figure 8 negative Pauli Fibonacci numbers can be constructed (see figure 18).

## 9. Jacobsthal Numbers and Pauli Jacobsthal Numbers

If every pair of rabbits produces two pairs of rabbits of the next generation, Jacobsthal numbers [13, p. 951] will emerge (see figure 19).

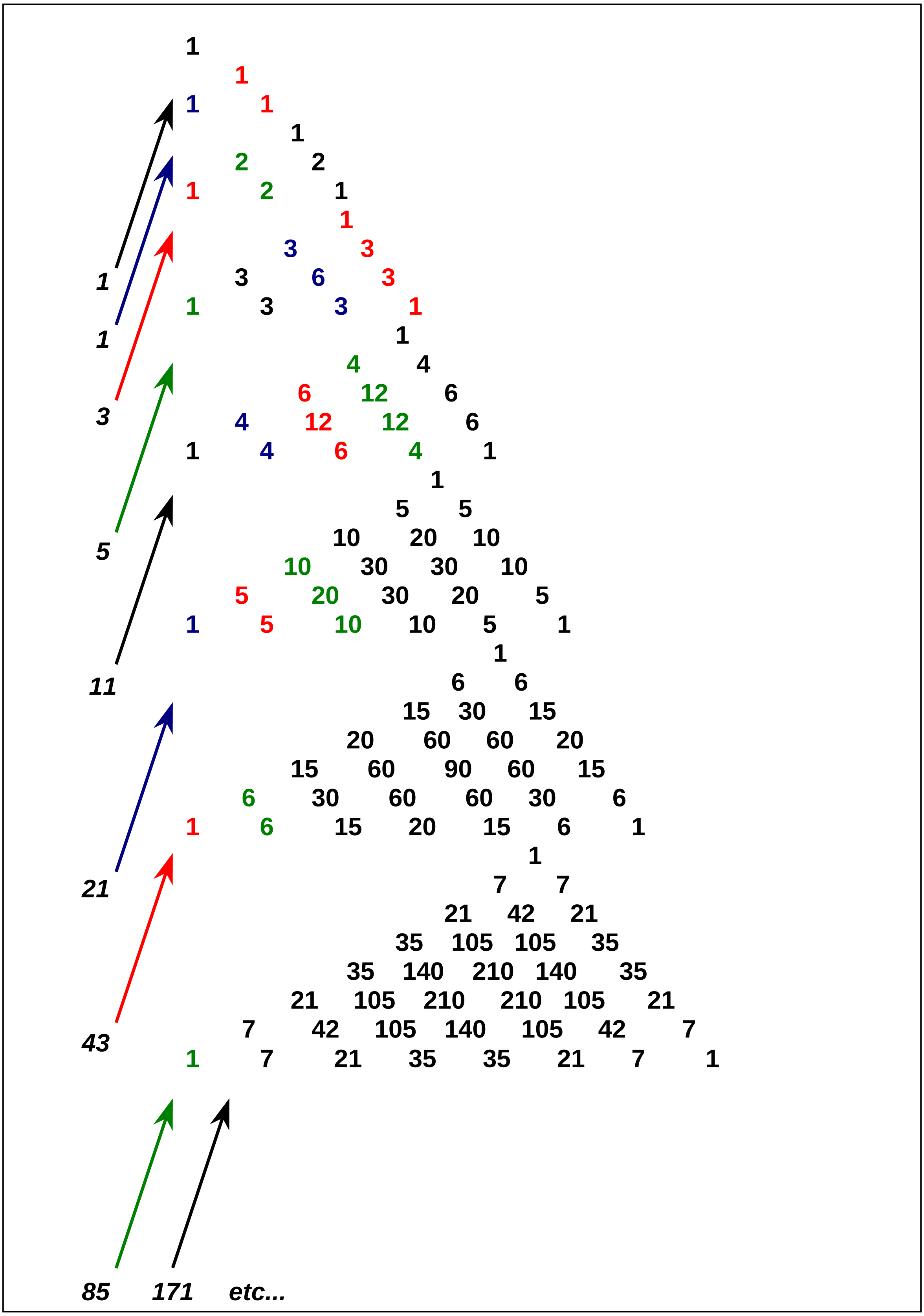

Fig. 19: Construction of positive Jacobsthal numbers in three-dimensional Pascal space.

They satisfy the recurrence relation

$$J_n = J_{n-1} + 2\, J_{n-2} \tag{28}$$

According to the procedure in [12, p. 18 (1.18)] Binet's formula for Jacobsthal numbers can be found by solving the characteristic equation

$$2 + q = q^2 \tag{29}$$

which gives

$$q_{1/2} = \frac{1 \pm \sqrt{1+8}}{2} = \frac{1 \pm 3}{2} \tag{30}$$

Therefore the closed form for Jacobsthal numbers is

$$J_n = \frac{1}{3}\left[\left(\frac{1+3}{2}\right)^n - \left(\frac{1-3}{2}\right)^n\right] = \frac{1}{3}\left(2^n - (-1)^n\right) \tag{31}$$

With the recurrence relation $J_{n-2} = (J_n - J_{n-1})/2$ or directly with Binet's formula (31) negative Jacobsthal numbers can be calculated:

$$\begin{aligned} J_3 &= 3 \\ J_2 &= 1 \\ J_1 &= 1 \\ J_0 &= 0 \\ J_{-1} &= \frac{1}{2} \\ J_{-2} &= -\frac{1}{4} \\ J_{-3} &= \frac{3}{8} \\ J_{-4} &= -\frac{5}{16} \\ J_{-5} &= \frac{11}{32} \\ J_{-6} &= -\frac{21}{64} \end{aligned} \tag{32}$$

But where are these numbers in figure 4? In analogy to the construction of Fibonacci numbers (see figure 14), Jacobsthal numbers should be found (see figure 20), but unfortunately the series don't converge. Therefore some imagination is needed to reconstruct negative Jacobsthal numbers.

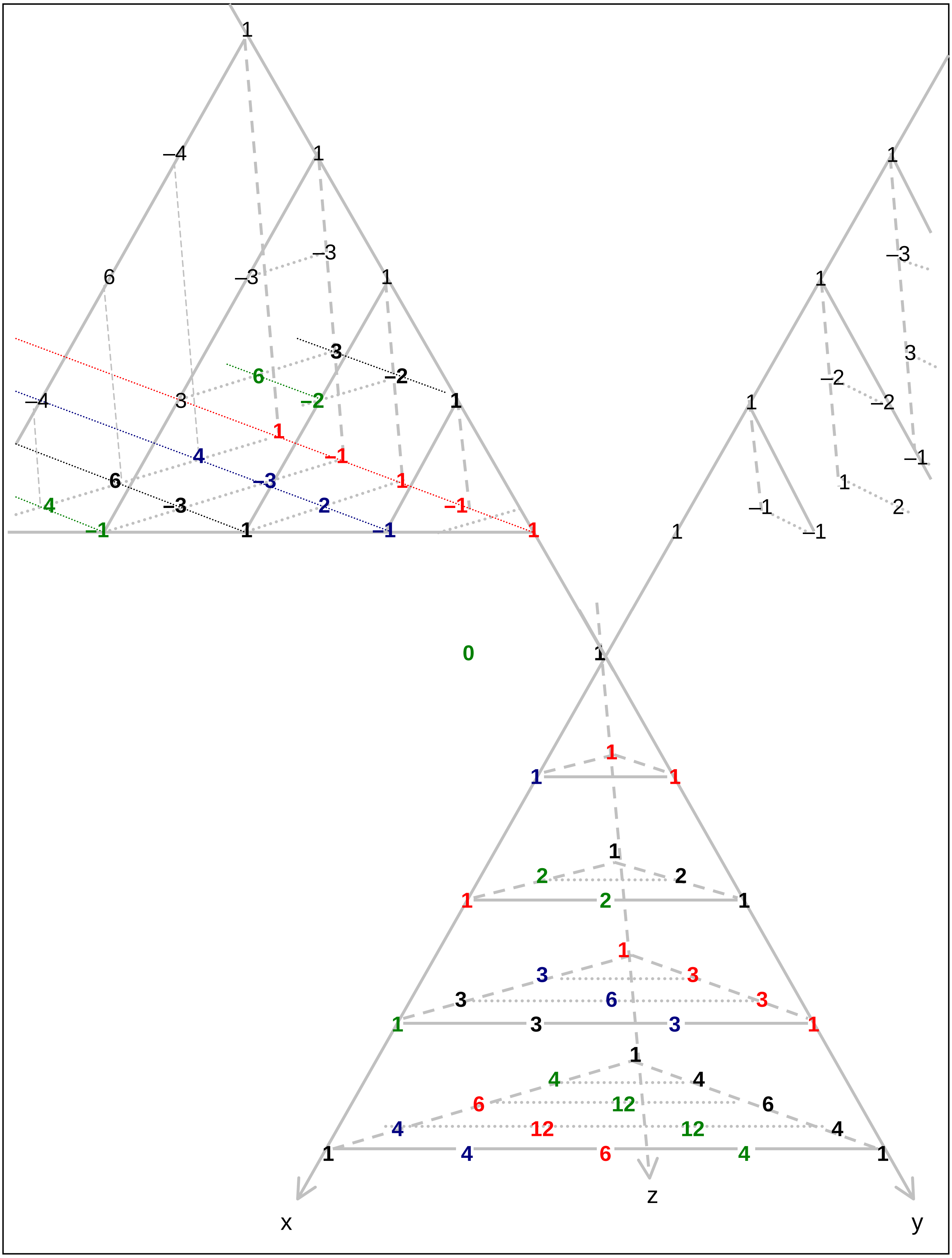

Fig. 20: Some imagination is needed to reconstruct negative Jacobsthal numbers.

Euler had this imagination, as Knopp wrote: „Und so ist es nicht zu verwundern, dass z.B. Euler die geometrische Reihe

$$1+x+x^2+\ldots=\frac{1}{1-x}$$

auch für $x=-1$ oder $x=-2$ noch gelten lässt und also unbedenklich

$$1-1+1-1+-\ldots=\frac{1}{2}$$

oder

$$1-2+2^2-2^3+-\ldots=\frac{1}{3}$$

setzt, und entsprechend etwa aus $\left(\frac{1}{1-x}\right)^2=1+2x+3x^2+\ldots$ die Gleichung

$$1-2+3-4+-\ldots=\frac{1}{4}$$

herleitet und vieles mehr." [9, p. 474] More reasons, why

$$\begin{aligned} J_{-1} &= 1-1+1-1+-\ldots=\frac{1}{2} \\ J_{-2} &= -1+2-3+4-+\ldots=-\frac{1}{4} \\ J_{-3} &= (1-3+6-12+-\ldots)+(1-2+3-4+-\ldots)=\frac{1}{8}+\frac{1}{4}=\frac{3}{8} \\ &\text{etc}\ldots \end{aligned} \tag{33}$$

can be found in [9, p. 473-510].

Pauli Jacobsthal numbers can be constructed in a similar way. If every pair of rabbits produces two pairs of rabbits of the next generation, and every second generation of parents develop into killer rabbits, Pauli Jacobsthal numbers will emerge (see figure 21 and figure 22).

Like in relation (24) Pauli Jacobsthal numbers $PJ_n$ show a twofold symmetry

$PJ_{-3} = J_0 = 0$ $PJ_{-1} = J_1 = 1$ $PJ_1 = J_2 = 1$ $PJ_3 = J_3 = 3$ $PF_5 = J_4 = 5$ $PJ_7 = J_5 = 11$ $PJ_9 = J_6 = 21$

$PJ_{-2} = J_{-1} = 1/2$ $PJ_0 = J_0 = 0$ $PJ_2 = J_1 = 1$ $PJ_4 = J_2 = 1$ $PJ_6 = J_3 = 3$ $PJ_8 = J_4 = 5$ (34)

even in the negative region if something like

$$\begin{aligned} PJ_{-1} &= 1+1-1-1+1+1-1-1++--\ldots \\ &= (1-1+1-1+-\ldots)+(1-1+1-1+-\ldots)+=\frac{1}{2}+\frac{1}{2}=1 \end{aligned} \tag{35}$$

or

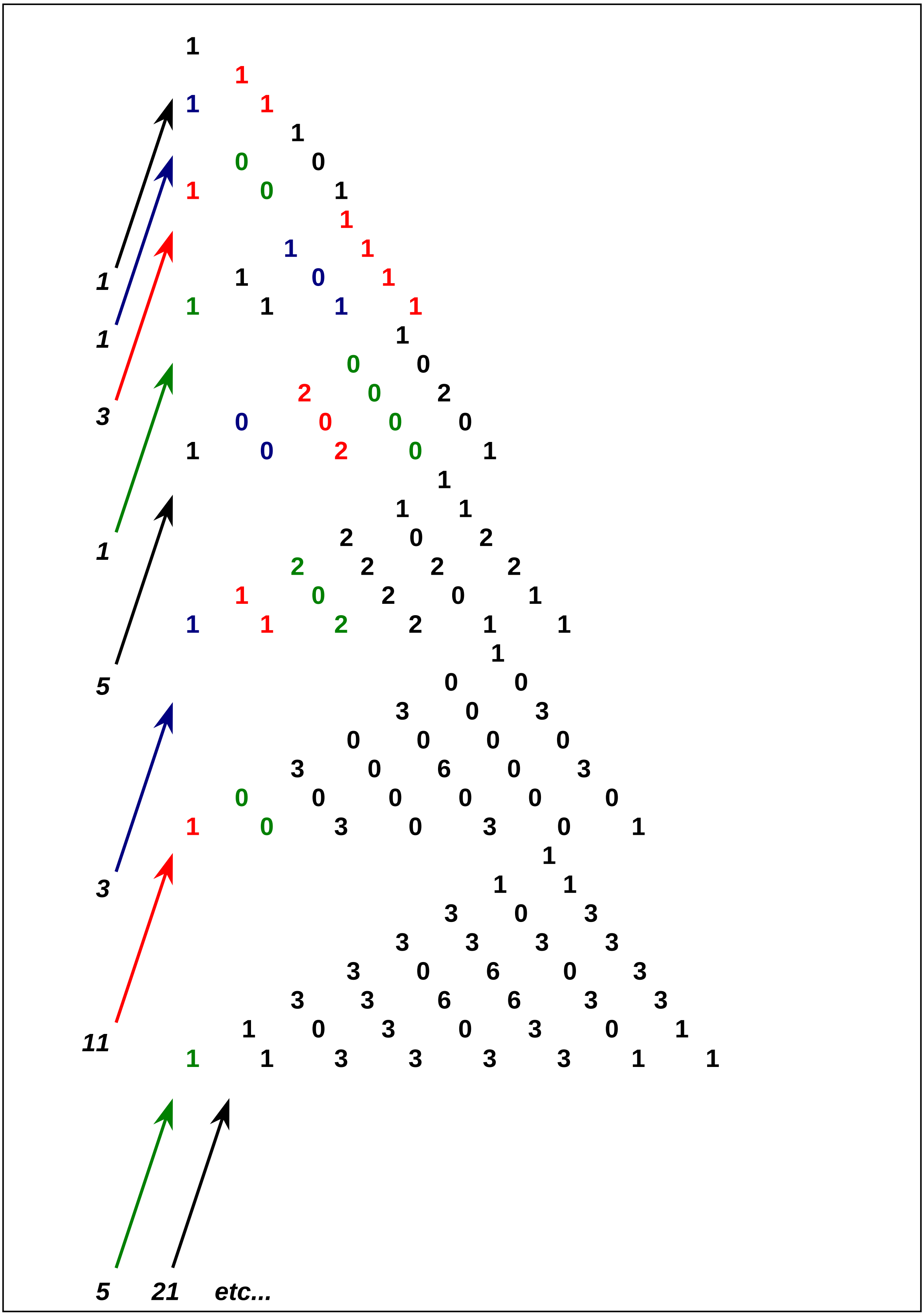

Fig. 21: Construction of positive Pauli Jacobsthal numbers in three-dimensional Pauli Pascal space.

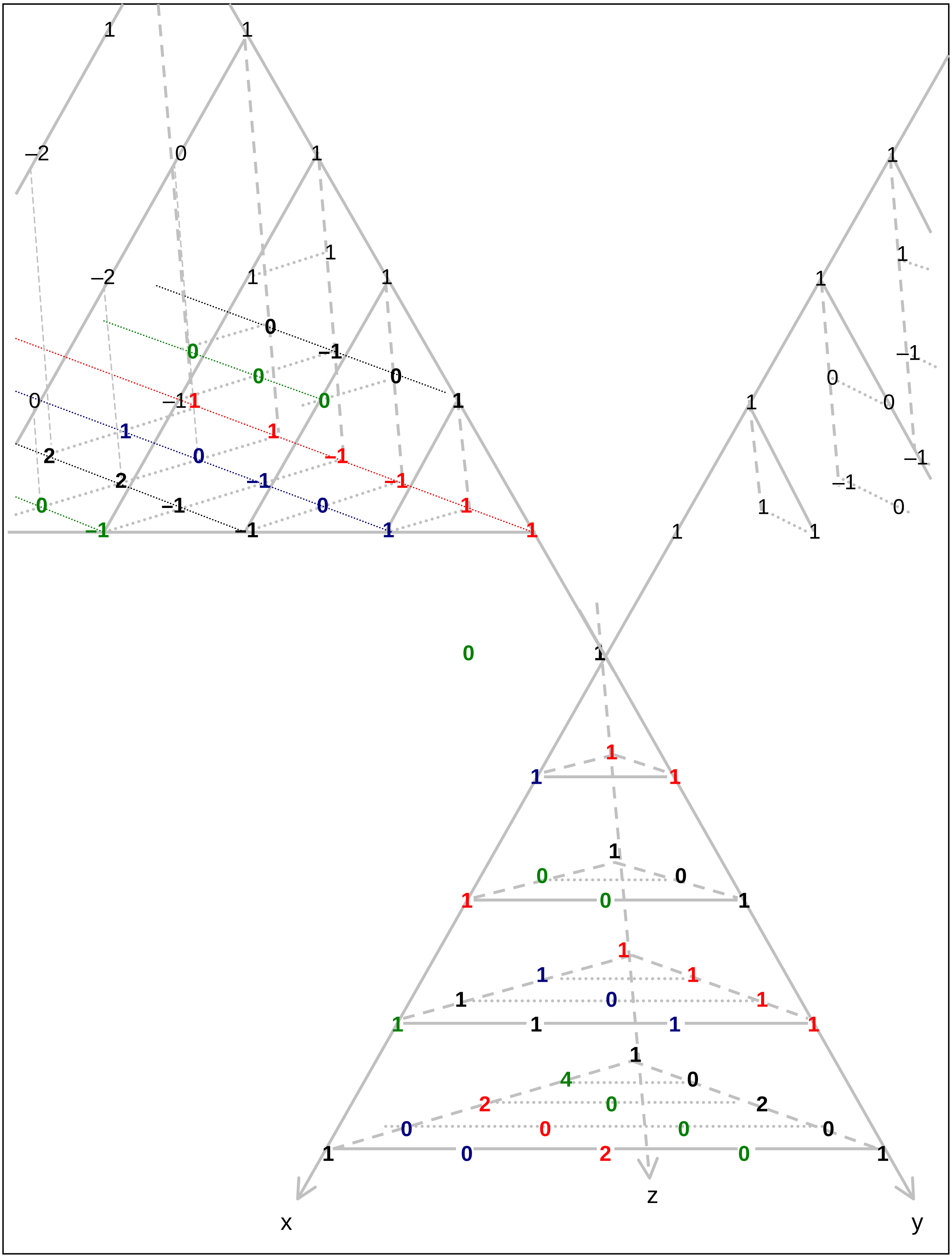

Fig. 22: Negative Pauli Jacobsthal numbers are supposed to be on the top left corner.

$$PJ_{-3} = (-1-1+2+2--++...)+(1-1+-...)$$
$$= (-1+2-...)+(-1+2-+...)+(1-1+-...) = -\frac{1}{4}-\frac{1}{4}+\frac{1}{2} = 0 \qquad (36)$$

seems to be valid.

In closed form Pauli Jacobsthal numbers are

$$PJ_n = \frac{1}{6}\left(\sqrt{2}^n + \sqrt{2}^{n+3} + \left(-\sqrt{2}\right)^n + \left(-\sqrt{2}\right)^{n+3} - i^n - i^{n+3} - (-i)^n - (-i)^{n+3}\right)$$
$$= \frac{1}{6}\left[\left(1+\sqrt{8}\right)\sqrt{2}^n + \left(1-\sqrt{8}\right)\left(-\sqrt{2}\right)^n - (1-i)i^n - (1+i)(-i)^n\right] \qquad (37)$$

## 10. Fibonacci Numbers and Pauli Fibonacci Numbers of Higher Order

Every pair of rabbits produces n pairs of rabbits of the next generation. If n is called the order of Fibonacci numbers, ordinary Fibonacci numbers are Fibonacci numbers of order one. Jacobsthal numbers are Fibonacci numbers of order two.

Of course it is possible that fast-growing pairs of rabbits produce more than two pairs of rabbits of the next generation simultaneously. Then Fibonacci numbers of higher orders can be found. They can be constructed by adding appropriate numbers of Fibonacci hyperpyramids resp. of the appropriate multinomial coefficients.

For example in four-dimensional Pascal space there is

$$\binom{n}{m_1 \quad m_2 \quad m_3 \quad n-m_1-m_2-m_3} = \frac{n!}{m_1!\cdot m_2!\cdot m_3!\cdot (n-m_1-m_2-m_3)!} \qquad (38)$$

and therefore Fibonacci numbers of order three are:

$$F_1^{(3)} = \binom{0}{0 \quad 0 \quad 0 \quad 0} = 1$$

$$F_2^{(3)} = \binom{1}{1 \quad 0 \quad 0 \quad 0} = 1$$

$$F_3^{(3)} = \binom{2}{2 \quad 0 \quad 0 \quad 0} + \binom{1}{0 \quad 1 \quad 0 \quad 0} + \binom{1}{0 \quad 0 \quad 1 \quad 0} + \binom{1}{0 \quad 0 \quad 0 \quad 1}$$
$$= 1+1+1+1 = 4$$

$$F_4^{(3)} = \binom{3}{3 \quad 0 \quad 0 \quad 0} + \binom{2}{1 \quad 1 \quad 0 \quad 0} + \binom{2}{1 \quad 0 \quad 1 \quad 0} + \binom{2}{1 \quad 0 \quad 0 \quad 1} \qquad (39)$$
$$= 1+2+2+2 = 7$$

$$
\begin{aligned}
F_5^{(3)} &= \binom{4}{4\;0\;0\;0}+\binom{3}{2\;1\;0\;0}+\binom{3}{2\;0\;1\;0}+\binom{3}{2\;0\;0\;1}+\binom{2}{0\;2\;0\;0}\\
&\quad+\binom{2}{0\;0\;2\;0}+\binom{2}{0\;0\;0\;2}+\binom{2}{0\;1\;1\;0}+\binom{2}{0\;1\;0\;1}+\binom{2}{0\;0\;1\;1}\\
&= 1+3+3+3+1+1+1+2+2+2=19
\end{aligned}
$$

$$
\begin{aligned}
F_6^{(3)} &= \binom{5}{5\;0\;0\;0}+\binom{4}{3\;1\;0\;0}+\binom{4}{3\;0\;1\;0}+\binom{4}{3\;0\;0\;1}+\binom{3}{1\;2\;0\;0}\\
&\quad+\binom{3}{1\;0\;2\;0}+\binom{3}{1\;0\;0\;2}+\binom{3}{1\;1\;1\;0}+\binom{3}{1\;1\;0\;1}+\binom{3}{1\;0\;1\;1}\\
&= 1+4+4+4+3+3+3+6+6+6=40
\end{aligned}
$$

$$
\begin{aligned}
F_7^{(3)} &= \binom{6}{6\;0\;0\;0}+\binom{5}{4\;1\;0\;0}+\binom{5}{4\;0\;1\;0}+\binom{5}{4\;0\;0\;1}+\binom{4}{2\;2\;0\;0}\\
&\quad+\binom{4}{2\;0\;2\;0}+\binom{4}{2\;0\;0\;2}+\binom{4}{2\;1\;1\;0}+\binom{4}{2\;1\;0\;1}+\binom{4}{2\;0\;1\;1}\\
&\quad+\binom{3}{0\;3\;0\;0}+\binom{3}{0\;0\;3\;0}+\binom{3}{0\;0\;0\;3}+\binom{3}{0\;2\;1\;0}+\binom{3}{0\;2\;0\;1}\\
&\quad+\binom{3}{0\;1\;2\;0}+\binom{3}{0\;0\;2\;1}+\binom{3}{0\;1\;0\;2}+\binom{3}{0\;0\;1\;2}+\binom{3}{0\;1\;1\;1}\\
&= 1+5+5+5+6+6+6+12+12+12+1+1+1+3+3+3+3+3+3+6=97
\end{aligned}
$$

$$
\begin{aligned}
F_8^{(3)} &= \binom{7}{7\;0\;0\;0}+\binom{6}{5\;1\;0\;0}+\binom{6}{5\;0\;1\;0}+\binom{6}{5\;0\;0\;1}+\binom{5}{3\;2\;0\;0}\\
&\quad+\binom{5}{3\;0\;2\;0}+\binom{5}{3\;0\;0\;2}+\binom{5}{3\;1\;1\;0}+\binom{5}{3\;1\;0\;1}+\binom{5}{3\;0\;1\;1}\\
&\quad+\binom{4}{1\;3\;0\;0}+\binom{4}{1\;0\;3\;0}+\binom{4}{1\;0\;0\;3}+\binom{4}{1\;2\;1\;0}+\binom{4}{1\;2\;0\;1}\\
&\quad+\binom{4}{1\;1\;2\;0}+\binom{4}{1\;0\;2\;1}+\binom{4}{1\;1\;0\;2}+\binom{4}{1\;0\;1\;2}+\binom{4}{1\;1\;1\;1}\\
&= 1+6+6+6+10+10+10+20+20+20+4+4+4+12+12+12+12+12+12+24\\
&= 217
\end{aligned}
$$

etc.

These Fibonacci numbers of third order satisfy the recurrence relation

$$F_{n-1}^{(3)} + 3F_{n-2}^{(3)} = F_n^{(3)} \tag{40}$$

According to the procedure in [12, p. 18 (1.18)] Binet's formula for these numbers can be found by solving the characteristic equation

$$3 + q = q^2 \tag{41}$$

which gives

$$q_{1/2} = \frac{1 \pm \sqrt{1+12}}{2} = \frac{1 \pm \sqrt{13}}{2} \tag{42}$$

Therefore the closed form for Fibonacci numbers of third order is

$$F_n^{(3)} = \frac{1}{\sqrt{13}}\left[\left(\frac{1+\sqrt{13}}{2}\right)^n - \left(\frac{1-\sqrt{13}}{2}\right)^n\right] \tag{43}$$

This recurrence relation is again valid for negative n:

$$\begin{aligned}
F_{-n}^{(3)} &= \frac{1}{\sqrt{13}}\left[\left(\frac{1+\sqrt{13}}{2}\right)^{-n} - \left(\frac{1-\sqrt{13}}{2}\right)^{-n}\right] \\
&= \frac{1}{\sqrt{13}}\left[\left(\frac{2}{1+\sqrt{13}}\right)^n - \left(\frac{2}{1-\sqrt{13}}\right)^n\right] \\
&= \frac{1}{\sqrt{13}}\left[\left(\frac{2\left(1-\sqrt{13}\right)}{1-13}\right)^n - \left(\frac{2\left(1+\sqrt{13}\right)}{1-13}\right)^n\right] \\
&= \frac{1}{\sqrt{13}}\left[\left(-\frac{1-\sqrt{13}}{6}\right)^n - \left(-\frac{1+\sqrt{13}}{6}\right)^n\right] \\
&= \frac{1}{\sqrt{13}\cdot(-3)^n}\left[\left(\frac{1+\sqrt{13}}{2}\right)^n - \left(\frac{1-\sqrt{13}}{2}\right)^n\right] = \frac{1}{(-3)^n}\cdot F_n^{(3)}
\end{aligned} \tag{44}$$

These negative Fibonacci numbers of third order can be found, if multinomial coefficients are calculated with the help of the Gamma function

$$\binom{n}{m_1 \quad m_2 \quad m_3 \quad n-m_1-m_2-m_3} = \lim_{h\to 0}\frac{\Gamma(n+1+h)}{\Gamma(m_1+1+h)\cdot\Gamma(m_2+1)\cdot\Gamma(m_3+1)\cdot\Gamma(n-m_1-m_2-m_3+1)}$$

(for $n, m_1 < 0$) (45)

and using Eulerian methods [9, p. 473-510]. Then

$$F_0^{(3)} = \binom{0}{1 \quad -1 \quad 0 \quad 0} = 0$$

$$
\begin{aligned}
F_{-1}^{(3)} &= \begin{pmatrix} & -1 & & \\ 0 & -1 & 0 & 0 \end{pmatrix} + \begin{pmatrix} & -1 & & \\ 0 & -2 & 1 & 0 \end{pmatrix} + \begin{pmatrix} & -1 & & \\ 0 & -2 & 0 & 1 \end{pmatrix} \\
&\quad + \begin{pmatrix} & -1 & & \\ 0 & -3 & 2 & 0 \end{pmatrix} + \begin{pmatrix} & -1 & & \\ 0 & -3 & 1 & 1 \end{pmatrix} + \begin{pmatrix} & -1 & & \\ 0 & -3 & 0 & 2 \end{pmatrix} \\
&\quad + \begin{pmatrix} & -1 & & \\ 0 & -4 & 3 & 0 \end{pmatrix} + \begin{pmatrix} & -1 & & \\ 0 & -4 & 2 & 1 \end{pmatrix} + \begin{pmatrix} & -1 & & \\ 0 & -4 & 1 & 2 \end{pmatrix} + \begin{pmatrix} & -1 & & \\ 0 & -4 & 0 & 3 \end{pmatrix} + \\
&\quad + \begin{pmatrix} & -1 & & \\ 0 & -5 & 4 & 0 \end{pmatrix} + \begin{pmatrix} & -1 & & \\ 0 & -5 & 3 & 1 \end{pmatrix} + \begin{pmatrix} & -1 & & \\ 0 & -5 & 2 & 2 \end{pmatrix} + \begin{pmatrix} & -1 & & \\ 0 & -5 & 1 & 3 \end{pmatrix} + \begin{pmatrix} & -1 & & \\ 0 & -5 & 0 & 4 \end{pmatrix} + \ldots \\
&= 1 - 1 - 1 + 1 + 2 + 1 - 1 - 3 - 3 - 1 + 1 + 4 + 6 + 4 + 1 - + \ldots \\
&= 1 - 2 + 4 - 8 + 16 - + \ldots \\
&= 1 - 2^1 + 2^2 - 2^3 + 4^4 - + \ldots \\
&= \frac{1}{3} \qquad (46)
\end{aligned}
$$

$$F_{-2}^{(3)} = -\frac{1}{9}$$

$$F_{-3}^{(3)} = \frac{4}{27}$$

$$F_{-4}^{(3)} = -\frac{7}{81}$$

$$F_{-5}^{(3)} = \frac{19}{243}$$

$$F_{-6}^{(3)} = -\frac{40}{729}$$

$$F_{-7}^{(3)} = \frac{97}{2187}$$

In a similar way Pauli Fibonacci numbers of third order can be found. In closed form they are:

$$
\begin{aligned}
PF_n^{(3)} = \frac{1}{2\sqrt{13}} \Bigg[ & \sqrt{\frac{1+\sqrt{13}}{2}}^{\,n} + \left(-\sqrt{\frac{1+\sqrt{13}}{2}}\right)^{n} - \sqrt{\frac{1+\sqrt{13}}{2}}^{\,n} - \left(-\sqrt{\frac{1+\sqrt{13}}{2}}\right)^{n} \\
& + \sqrt{\frac{1+\sqrt{13}}{2}}^{\,n+3} + \left(-\sqrt{\frac{1+\sqrt{13}}{2}}\right)^{n+3} - \sqrt{\frac{1+\sqrt{13}}{2}}^{\,n+3} - \left(-\sqrt{\frac{1+\sqrt{13}}{2}}\right)^{n+3} \Bigg] \qquad (47)
\end{aligned}
$$

Fibonacci numbers of order k satisfy the recurrence relation

$$F_{n-1}^{(k)} + kF_{n-2}^{(k)} = F_n^{(k)} \qquad (48)$$

According to the procedure in [12, p. 18 (1.18)] Binet's formula for these numbers can be found by solving

$$k + q = q^2 \qquad (49)$$

which gives

$$q_{1/2} = \frac{1 \pm \sqrt{1+4k}}{2} \tag{50}$$

Therefore the closed form for Fibonacci numbers of order k is

$$F_n^{(k)} = \frac{1}{\sqrt{1+4k}} \left[ \left( \frac{1+\sqrt{1+4k}}{2} \right)^n - \left( \frac{1-\sqrt{1+4k}}{2} \right)^n \right] \tag{51}$$

Friendly squares with odd integers will result when

$$\begin{aligned}
k_0 = 0 &\quad \Rightarrow \quad q_{1/2} = \frac{1 \pm \sqrt{1+0}}{2} = \frac{1}{2}(1 \pm 1) \\
k_2 = 2 &\quad \Rightarrow \quad q_{1/2} = \frac{1 \pm \sqrt{1+8}}{2} = \frac{1}{2}(1 \pm 3) \\
k_4 = 6 &\quad \Rightarrow \quad q_{1/2} = \frac{1 \pm \sqrt{1+24}}{2} = \frac{1}{2}(1 \pm 5) \\
k_6 = 12 &\quad \Rightarrow \quad q_{1/2} = \frac{1 \pm \sqrt{1+48}}{2} = \frac{1}{2}(1 \pm 7) \\
k_8 = 20 &\quad \Rightarrow \quad q_{1/2} = \frac{1 \pm \sqrt{1+80}}{2} = \frac{1}{2}(1 \pm 9) \\
k_{10} = 30 &\quad \Rightarrow \quad q_{1/2} = \frac{1 \pm \sqrt{1+120}}{2} = \frac{1}{2}(1 \pm 11) \\
k_{m-1} = \frac{m+1}{2} \cdot \frac{m-1}{2} &\quad \Rightarrow \quad q_{1/2} = \frac{1 \pm \sqrt{1+(m^2-1)}}{2} = \frac{1}{2}(1 \pm m)
\end{aligned} \tag{52}$$

In a similar way Pauli Fibonacci numbers of higher order can be constructed.

## 11. Outlook

In theoretical physics we can never be sure whether the reality we find outside us is something which comes directly from nature or whether it is something we have invented when inventing mathematics. As the footnotes indicate, there might be some physical consequences of Taylor expanding formulae in physics. We all live in Pascal space because human[3] theoretical physicists use the Taylor expansion.

And there are interesting mathematical consequences as Ramanujan's $_1\psi_1$–summation formula may be decoded in Pascal and Pauli Pascal space too.

---

[3] "If there were extraterrestrial aliens, they very probably would invent totally different mathematics," says Kovalevskaya-prize-winning Matilde Marcolli. („Wenn es außerirdische Lebewesen gäbe, dann würden sie höchstwahrscheinlich auch eine vollkommen andere Mathematik erfinden." [11, p. 80])